\documentclass[onecolumn,final,a4page]{siamltex}

\usepackage{amsmath}       
\usepackage{amssymb}       
\usepackage{chapterbib}    
\usepackage{epsfig}        
\usepackage{color}
\usepackage{comment}
\usepackage{bm,bbm}
\usepackage{overpic}

\def\Frac{\displaystyle \frac}
\def\Int{ \displaystyle \int}
\def\Sum{ \displaystyle \sum}

\def\div{{ \rm div}}

\def\DIV{{\cal DIV}}
\def\GRAD{{\cal GRAD}}

\def\O{\Omega}

\def\bn{{\bf n}}

\def\br{{\bf r}}
\def\bu{{\bf u}}
\def\bv{{\bf v}}
\def\bw{{\bf w}}
\def\bx{{\bf x}}

\def\cF{{\cal F}}

\def\cP{{\cal P}}

\def\cS{{\cal S}}

\renewcommand{\theequation}{\thesection.\arabic{equation}}

\newtheorem{assumption}{Assumption}[section]
\newtheorem{remark}{Remark}[section]

\def\k{k}

\usepackage[normalem]{ulem}
\definecolor{violet}{rgb}{0.580,0.,0.827}

\usepackage{textpos}
\usepackage{cite}

\title{Convergence analysis of the mimetic finite difference method for elliptic problems with
  staggered discretizations of diffusion coefficients}

\author{
  G. Manzini   \footnotemark[1],\quad
  K. Lipnikov  \footnotemark[1],\quad
  J. D. Moulton\footnotemark[1],\quad
  M. Shashkov  \footnotemark[2]
}
  
\begin{document}

\maketitle

\pagestyle{myheadings}
\thispagestyle{plain}

\markboth
{ G.~Manzini, L.~Lipnikov, J.~D.~Moulton, M.~Shashkov}
{The MFD method for elliptic problems with staggered discretizations of diffusion coefficients}

\renewcommand{\thefootnote}{\fnsymbol{footnote}}

\footnotetext[1]{Applied Mathematics and Plasma Physics Group, Theoretical Division,
  Los Alamos National Laboratory, \{lipnikov,gmanzini,moulton\}@lanl.gov}

\footnotetext[2]{XCP-4 Group, Computational Physics Division, 
  Los Alamos National Laboratory, shashkov@lanl.gov}

\renewcommand{\thefootnote}{\arabic{footnote}}

\markboth
{G.~Manzini, L.~Lipnikov, J.~D.~Moulton, M.~Shashkov}
{The MFD method for elliptic problems with staggered discretizations of diffusion coefficients}

\begin{abstract}
  We study the convergence of the new family of mimetic finite
  difference schemes for linear diffusion problems recently proposed
  in~\cite{Lipnikov-Manzini-Moulton-Shashkov:2016}.
  In contrast to the conventional approach, the diffusion coefficient
  enters both the primary mimetic operator, i.e., the discrete
  divergence, and the inner product in the space of gradients.
  The diffusion coefficient is therefore evaluated on different mesh
  locations, i.e., inside mesh cells and on mesh faces.
  Such a staggered discretization may provide the flexibility
  necessary for future development of efficient numerical schemes for
  nonlinear problems, especially for problems with degenerate
  coefficients.
  These new mimetic schemes preserve symmetry and
  positive-definiteness of the continuum problem, which allow us to
  use efficient algebraic solvers such as the preconditioned Conjugate
  Gradient method.
  We show that these schemes are inf-sup stable and establish a priori
  error estimates for the approximation of the scalar and vector
  solution fields.
  Numerical examples confirm the convergence analysis and the
  effectiveness of the method in providing accurate approximations.
\end{abstract}




\raggedbottom

\newcommand{\RED} [1]{{\color{red}#1}}
\newcommand{\BLUE}[1]{{\color{blue}#1}}
\newcommand{\MAGENTA}[1]{{\color{magenta}#1}}

\newcommand{\SchemeI} {\textsf{Trace}{}}
\newcommand{\SchemeII}{\textsf{Upwind}{}}

\newcommand{\hT}{h_T}
\newcommand{\rT}{r_T}
\newcommand{\aStar}{a^{\star}}
\newcommand{\astar}{a_{\star}}
\newcommand{\bstar}{b_{\star}}

\newcommand{\Ns}{{\cal N}_s}

\newcommand{\bvRT}{\bv_{RT_0}}
\newcommand{\CRTz}{C_{RT_0}}

\newcommand{\bwh}{\bw_h}
\newcommand{\barkc}{\bar{\k}^c}

\newcommand{\bilAh} [2]{a_{h}(#1,#2)}
\newcommand{\bilAhc}[2]{a_{h,c}(#1,#2)}
\newcommand{\bilA}  [2]{a(#1,#2)}
\newcommand{\bilAc} [2]{a_{c}(#1,#2)}

\newcommand{\STAB}  [2]{\textsf{stab}(#1,#2)}
\newcommand{\STABc} [2]{\textsf{stab}_c(#1,#2)}

\newcommand{\Pizc}{\Pi^0_{c}}

\newcommand{\DIVk} {\DIV^{\k}}
\newcommand{\DIVkc}{\DIV^{\k}_{c}}

\newcommand{\CIp}{C^*_{Ip}}
\newcommand{\CAg}{C^*_{Ag}}

\newcommand{\II}{I\hspace{-0.35mm}I}

\newcommand{\Hdiv}{H_{\div}}
\newcommand{\cf}[1]{{#1}^c_f}
\newcommand{\cof}[1]{{#1}^{c_1}_f}
\newcommand{\ctf}[1]{{#1}^{c_2}_f}

\newcommand{\pc}{\partial c}

\newcommand{\CNST}[1]{C^*_{#1}}

\def\trait #1 #2 #3 {\vrule width #1pt height #2pt depth #3pt}
\def\fin{\hfill
        \trait .3 5 0
        \trait 5 .3 0
        \kern-5pt
        \trait 5 5 -4.7
        \trait 0.3 5 0
\medskip}
\newcommand{\ENDTOPIC}{\fin}
\newcommand{\ENDPROOF}{\fin}
\newcommand{\BEGINPROOF}{\emph{Proof}.~~}

\newcommand{\FTN}[1]{\footnote{\textbf{#1}}}
\newcommand{\ADD}[1]{\textbf{#1}}
\newcommand{\BAD}{\textbf{[This is bad!!!]}}
\newcommand{\EOD}{\end{document}} 

\newcommand{\scf}{\sigma^c_f}
\newcommand{\scof}{\sigma^{c_1}_f}
\newcommand{\sctf}{\sigma^{c_2}_f}
\newcommand{\scfo}{\sigma^c_{f_1}}
\newcommand{\scfn}{\sigma^c_{f_n}}

\newcommand{\matKc}{\widetilde{{\cal K}}_{c}}
\newcommand{\kc}  {\k^c}
\newcommand{\ktcf}{\widetilde{\k}^c_f}
\newcommand{\kcf} {\k^c_f}
\newcommand{\kci} {\k^{c_i}}
\newcommand{\kcof}{\k^{c_1}_f}
\newcommand{\kctf}{\k^{c_2}_f}
\newcommand{\kcif}{\k^{c_i}_f}
\newcommand{\kk}  {\k}

\newcommand{\ktcof}{\widetilde{\k}^{c_1}_f}
\newcommand{\ktctf}{\widetilde{\k}^{c_2}_f}
\newcommand{\ktcif}{\widetilde{\k}^{c_i}_f}

\newcommand{\kbf} {\overline{\k}_f}

\newcommand{\ucf}{u^c_f}
\newcommand{\vcf}{v^c_f}
\newcommand{\wcf}{w^c_f}

\newcommand{\ucof}{u^{c_1}_f}
\newcommand{\uctf}{u^{c_2}_f}

\newcommand{\vcof}{v^{c_1}_f}
\newcommand{\vctf}{v^{c_2}_f}

\newcommand{\buc}{\bu_c}
\newcommand{\bucI}{\bu_c^I}
 \newcommand{\bucII}{\bu_c^{\II}}
\newcommand{\buI} {\bu^I}
\newcommand{\buIc} {\bu^I_c}
\newcommand{\buII}{\bu^{\II}}
\newcommand{\buh}{\bu_h}

\newcommand{\buth}{\widetilde{\bu}_{h}}
\newcommand{\butc}{\widetilde{\bu}_{c}}

\newcommand{\bvc}{\bv_c}
\newcommand{\bvcI}{\bv_c^I}
\newcommand{\bvcII}{\bv_c^{\II}}
\newcommand{\bvI} {\bv^I}
\newcommand{\bvII}{\bv^{\II}}
\newcommand{\bvh}{\bv_h}
\newcommand{\bvq}{\bv_{\qh}}

\newcommand{\bvth}{\widetilde{\bv}_{h}}
\newcommand{\bvtc}{\widetilde{\bv}_{c}}

\newcommand{\bnf}{\bn_f}
\newcommand{\bncf}{\bn^c_f}

\newcommand{\brc}{\br_c}

\newcommand{\bvpsi}{\bv_{\psi}}
\newcommand{\psib}{\overline{\psi}}

\newcommand{\bvhp} {\bvh^{\psi}}
\newcommand{\bvthp}{\widetilde{\bv}_{h}^{\psi}}

\newcommand{\qh}{q_h}
\newcommand{\ph}{p_h}
\newcommand{\pI}{p^I}
\newcommand{\qI}{q^I}

\newcommand{\pIc}{p^I_c}
\newcommand{\qIc}{q^I_c}

\newcommand{\rcf}{r^c_f}

\newcommand{\hf}{h_f}
\newcommand{\hc}{h_c}
\newcommand{\Oh}{\O_h}

\newcommand{\dS}{dS}
\newcommand{\dV}{dV}

\newcommand{\Ph}{\cP_h}
\newcommand{\Fh}{\cF_h}

\newcommand{\cT}{\mathcal{T}}
\newcommand{\Th} {\cT_{h,c}}
\newcommand{\Thc}{\cT_h}

\newcommand{\Sh}{\cS_h}

\newcommand{\Fth}{\widetilde{\mathcal{F}}_h}

\newcommand{\Phc}{\cP_{h,c}}
\newcommand{\Fhc}{\cF_{h,c}}
\newcommand{\Shc}{\cT_{h,c}}

\newcommand{\scal}  [2]{\big(#1,#2\big)}
\newcommand{\bscal} [2]{\left<#1,#2\right>_{\Gamma}}
\newcommand{\bscalh}[2]{\left<#1,#2\right>_{\Gamma,h}}

\newcommand{\scalPh}[2]{\big[#1,#2\big]_{\Ph}}
\newcommand{\scalFh}[2]{\big[#1,#2\big]_{\Fh}}

\newcommand{\scalPhc}[2]{\big[#1,#2\big]_{\Ph,c}}
\newcommand{\scalFhc}[2]{\big[#1,#2\big]_{\Fh,c}}

\newcommand{\ScalPh}[2]{\Big[#1,#2\Big]_{\Ph}}
\newcommand{\ScalFh}[2]{\Big[#1,#2\Big]_{\Fh}}

\newcommand{\ScalPhc}[2]{\Big[#1,#2\Big]_{\Ph,c}}
\newcommand{\ScalFhc}[2]{\Big[#1,#2\Big]_{\Fh,c}}

\newcommand{\SCALPhc}[2]{\left[#1,#2\right]_{\Ph,c}}
\newcommand{\SCALFhc}[2]{\left[#1,#2\right]_{\Fh,c}}

\newcommand{\pone}   {p_1}
\newcommand{\qone}   {q_1}
\newcommand{\pzero}  {p_0}
\newcommand{\qzero}  {q_0}
\newcommand{\psione} {\psi_1}
\newcommand{\psizero}{\psi_0}
\newcommand{\bvzero} {\bv_0}
\newcommand{\buzero} {\bu_0}
\newcommand{\bvIzero}{\bv^I_0}

\newcommand{\nlen}{\hspace{-0.2mm}}
\newcommand{\hskp}{\hspace{0.2mm}}

\newcommand{\snorm}  [2]{\big|#1\big|_{#2}}
\newcommand{\Snorm}  [2]{\Big|#1\Big|_{#2}}
\newcommand{\SNORM}  [2]{|#1|_{#2}}

\newcommand{\norm}   [2]{\big|\nlen\big|#1\big|\nlen\big|_{#2}}
\newcommand{\Norm}   [2]{\Big|\nlen\nlen\Big|\hskp #1\hskp\Big|\nlen\nlen\Big|_{#2}}
\newcommand{\NORM}   [2]{|\nlen|#1|\nlen|_{#2}}

\newcommand{\Tnorm}  [2]{\big|\nlen\big|\nlen\big|#1\big|\nlen\big|\nlen\big|_{{#2}}}
\newcommand{\TNorm}  [2]{\Big|\nlen\nlen\Big|\nlen\nlen\Big|\hskp #1\hskp\Big|\nlen\nlen\Big|\nlen\nlen\Big|_{#2}}
\newcommand{\TNORM}  [2]{|\nlen|\nlen|#1|\nlen|\nlen|_{{}_{#2}}}

\newcommand{\abs}[1]{\left|#1\right|}
\newcommand{\ABS}[1]{\big|#1\big|}
\newcommand{\ASSUM}[1]{\textbf{#1}}

\newcommand{\PS}[1]{{\cal P}^{#1}}

\newcommand{\TEXTFONTA}[1]{\mathsf{#1}}
\newcommand{\matI} {\TEXTFONTA{I}}
\newcommand{\matN} {\TEXTFONTA{N}}
\newcommand{\matR} {\TEXTFONTA{R}}
\newcommand{\matM} {\TEXTFONTA{M}}
\newcommand{\matW} {\TEXTFONTA{W}}
\newcommand{\matG} {\TEXTFONTA{G}}
\newcommand{\matP} {\TEXTFONTA{P}}

\newcommand{\matNc}{\TEXTFONTA{N}_c}
\newcommand{\matRc}{\TEXTFONTA{R}_c}
\newcommand{\matWc}{\TEXTFONTA{W}_c}
\newcommand{\matGc}{\TEXTFONTA{G}_c}
\newcommand{\matPc}{\TEXTFONTA{P}_c}

\newcommand{\PGRAPH}[1]{\medskip\noindent\textbf{#1}}

\newcommand{\bepsh}{{\bm\varepsilon}_h}
\newcommand{\bepsc}{{\bm\varepsilon}_c}
\newcommand{\vepscf}{\varepsilon_f^c}
\newcommand{\vepscof}{\varepsilon^{c_1}_f}
\newcommand{\vepsctf}{\varepsilon^{c_2}_f}


\section{Introduction}

Complex geophysical subsurface and surface flows, including general
non-linear diffusion problems~\cite{Castor:2004} and moisture
transport in partially saturated porous media~\cite{Richards:1931} are
mathematically modeled through parabolic equations such as $\partial
\theta(p)/\partial t-\div\big(\k(p)\nabla p\big)=0$, where $\theta(p)$
and $k(p)$ are given nonlinear functions of the scalar unknown $p$.
The numerical approximation of this kind of equations is extremely
challenging when the diffusion coefficient $k(p)$ approaches zero due
to the non-linear dependence on $p$, or it presents very strong
discontinuities.
In such cases, the numerical approximation to $p$ becomes dramatically
inaccurate if $\k(p)$ is incorporated in the discrete form of the
equation through some kind of harmonic average of one-sided values of
$k^{-1}(p)$ at the mesh interfaces.
This fact is a major issue as it impacts almost all the discretization
methods in the literature that write the flux equation as
$\k^{-1}\bu=-\nabla p$.
This issue affects also the discretization of linear diffusion
problems where $\k$ is only a function of position, but may be
discontinuous or close to zero in some parts of the domain.
In the finite element (FE) and finite volume (FV) frameworks we
mention the mixed finite element
method~\cite{Boffi-Brezzi-Fortin:2013}, the `standard' mimetic finite
difference (MFD)
method~\cite{BeiraodaVeiga-Lipnikov-Manzini:2014,Lipnikov-Manzini-Shashkov:2014},
the gradient
scheme~\cite{Droniou-Eymard-Herbin:2016,Droniou-Nataraj-Devika:2016},
the hybrid and mixed finite volumes
method~\cite{Eymard-Gallouet-Herbin:2001,Droniou-Eymard:2006,Eymard-Gallouet-Herbin:2008,Eymard-Gallouet-Herbin:2000},
the hybrid high-order
method~\cite{DiPietro-Ern:2015,Cockburn-DiPietro-Ern:2016}, the mixed
weak Galerkin method~\cite{Wang-Ye:2013} and the mixed virtual element
method~\cite{Brezzi-Falk-Marini:2014,BeiraodaVeiga-Brezzi-Marini-Russo:2016}.
On the other hand, finite difference methods and finite volume methods
that approximate directly $\k\nabla p$ in the mass conservation
equation do not invert the diffusion coefficient and do not suffer of
this problem.
However, in these methods the symmetry of the discrete formulation is
typically lost and proving the coercivity, which implies that the
resulting matrix operator is positive definite, is a very hard and
sometimes impossible task~\cite{Eymard-Gallouet-Herbin:2000}.

Concerning numerical methods based on variational formulation, an
early success in addressing this issue and avoiding the inversion of
$\k(p)$ is found in~\cite{Arbogast-Wheeler-Yotov:1997,Arbogast:1998},
which proposed the expanded mixed FE method using two distinct vector
unknowns $\bu=-\nabla p$ and $\bv = \k\bu$.
However, this method has several drawbacks that motivate the current
work.
First, it is formulated only for finite element meshes of elements
with a few kind of geometric shapes, e.g., simplexes or quadrilaterals
in 2D and hexahedral and prismatic cells.
The current trend in the numerical treatment of partial differential
equations (PDEs) is toward applications using meshes with more general
polygonal and polyhedral elements.
The state of the art is reflected in the articles of the two recent
special
issues~\cite{Bellomo-Brezzi-Manzini:2014,BeiraodaVeiga-Ern:2016}.
Then, it employs only cell-centered diffusion coefficients but there
is strong evidence from practice that some sort of \emph{upwinding of
  the diffusion coefficient} is necessary for nonlinear problems.

For these reasons, in~\cite{Lipnikov-Manzini-Moulton-Shashkov:2016} we
proposed a new MFD formulation that is suitable to very general meshes
and uses a staggered representation of the diffusion coefficients at
the mesh interfaces.
In that first work, accuracy and robustness were assessed
experimentally for a set of steady-state linear diffusion problems and
a time-dependent parabolic problem with $\k$ approaching zero.
We emphasize that numerical and theoretical investigations on simpler
stationary linear problems are a necessary step for the proper design
of methods working on more complex time-dependent nonlinear problems.
In this paper, we support the numerical study of
~\cite{Lipnikov-Manzini-Moulton-Shashkov:2016} by theoretically
proving that our new MFD method is inf-sup stable, and, consequently,
well-posed, and is convergent when applied to the Poisson problem in
mixed form.
Convergence is proved by deriving first-order estimates for the scalar
and vector unknowns.
The extension of our methodology to time-dependent nonlinear problems
with degenerate coefficients will be the topic of future publications.

A mimetic method is specifically designed to preserve (or mimic)
essential mathematical and physical properties of the underlying PDEs
in the discrete setting.
For parabolic problems the essential properties may include the
corresponding conservation law, as well as the symmetry and
positive-definiteness of the underlying differential operator.
The MFD methodology both for the 
mixed formulation~\cite{%
BeiraodaVeiga-Droniou-Manzini:2011,%
BeiraodaVeiga-Gyrya-Lipnikov-Manzini:2009,%
BeiraodaVeiga-Lipnikov-Manzini:2009,%
BeiraodaVeiga-Lipnikov-Manzini:2010,%
BeiraodaVeiga-Manzini:2008b,%
BeiraodaVeiga-Manzini:2008,%
Brezzi-Lipnikov-Shashkov:2006,%
Brezzi-Lipnikov-Shashkov:2005,%
Brezzi-Lipnikov-Simoncini:2005,%
Cangiani-Manzini:2008,%
Cangiani-Manzini-Russo:2009,%
Cangiani-Gardini-Manzini:2011,%
Lipnikov-Morel-Shashkov:2004,%
Lipnikov-Manzini-Brezzi-Buffa:2011,%
Lipnikov-Manzini-Svyatskiy:2010} 
and the primal formulation~\cite{%
BeiraodaVeiga-Lipnikov-Manzini:2011,%
Brezzi-Buffa-Lipnikov:2009}%
~of elliptic problems has been the object of extensive development and
investigation during the last two decades, which proved its
effectiveness, accuracy and robustness.
For the interested readers, the main theoretical aspects in the
convergence analysis of the MFD method for elliptic PDEs are
summarized in the book~\cite{BeiraodaVeiga-Lipnikov-Manzini:2014}.
The book is complemented by two recently published review papers,
see~\cite{Lipnikov-Manzini-Shashkov:2014}
and~\cite{Gyrya-Lipnikov-Manzini-Svyatskiy:2014}.
In~\cite{Lipnikov-Manzini-Shashkov:2014} we review many known results
on Cartesian and curvilinear meshes for mathematical models that are
also non elliptic such as the Lagrangian hydrodynamics.
In~\cite{Gyrya-Lipnikov-Manzini-Svyatskiy:2014}, we review all known
optimization strategies that allows us to select schemes from the
mimetic family with superior properties, usually refereed to as the
\emph{mimetic optimization} or \emph{M-optimization}.
Such schemes may have a discrete maximum or minimum principle for
diffusion problems or show a significant reduction of the numerical
dispersion in wave propagation problems.

In the original mimetic framework, we discretize simultaneously pairs
of adjoint differential operators such as the divergence operator
\textsf{div}$(\cdot)$, and the flux operator $−\k\nabla(\cdot)$.
The divergence operator is chosen as the \emph{primary} operator and
is directly discretized consistently with the local Gauss divergence
theorem, while the discretization of the flux operator is
\emph{derived} from a discrete duality relation.
Instead, in the new MFD method the primary operator discretizes the
combined operator \textsf{div}$(\k\,\cdot\,)$ and the derived (dual)
gradient operator discretizes $\nabla(\cdot)$.
This alternative approach has two major consequences on the mimetic
formulation.
First, the mimetic inner product in the space of fluxes is weighted by
$\k$ instead of $\k^{-1}$ as in the original MFD method,
cf.~\cite{Brezzi-Lipnikov-Shashkov:2005,Brezzi-Lipnikov-Simoncini:2005}.
Second, a face-based representation of $\k$ is required in the
definition of the discrete divergence operator.
This staggered discretization allows us more freedom in the design of
a numerical method for vanishing or strongly discontinuous diffusion
coefficients $\k$ as we can use up to two distinct face values and
different ways to incorporate them in the discrete divergence
operator, e.g., through upwinding or arithmetic and harmonic
averaging.
It is also worth noting that the method resulting from this approach
in some specific cases includes other well-known ``classical''
schemes, e.g., the finite volume scheme using the two-point flux
approximation on orthogonal meshes, etc.
Finally, we note that our approach can be extended in a
straightforward way to the more general non-linear operator
$\div\big(\k(p)\mathbbm{K(\bx)}\nabla p\big)$ where $\mathbbm{K}(\bx)$
is a diffusion tensor dependent only on the position $\bx$ by
considering the splitting $\div(\k\,\cdot\,)$ and
$\mathbbm{K}\nabla(\cdot)$.
This generalization will be investigated in future works.

\medskip
The paper is organized as follows.
In Section~\ref{sec:mixed:formulation} we present the model problem.
In Section~\ref{sec:MFD:method} we formulate the new mimetic method
and discuss possible staggered approximations of the diffusion
coefficient at the mesh interfaces, e.g., first-order upwind and
arithmetic average of cell values of $\k$.
In Section~\ref{sec:error:estimates} we prove that the method is
well-posed and convergent and derive an a priori error estimate for
the approximation of the scalar and the gradient unknowns.
In Section~\ref{sec:numerics} we assess the behavior of the method
through numerical experiments.
In Section~\ref{sec:conclusions} we offer our final conclusions.

\subsection{Notation}
\medskip
Throughout the paper, we use the standard notation of Sobolev spaces,
cf.~\cite{Adams-Fournier:2003}.
In particular, let $\omega$ denote a domain in one or several
dimensions.
Then, $L^p(\omega)$, for any real $p$ such that $1\leq p<\infty$, is
the Sobolev space of $p$-integrable scalar functions and
$L^{\infty}(\omega)$ is the space of (essentially) bounded functions
defined on $\omega$; $W^{m,p}(\omega)$, for any integer $m\geq 1$ and
$1\leq p\leq\infty$, is the Sobolev space of functions in
$L^{p}(\omega)$ with all derivatives up to order $m$ also in
$L^p(\omega)$.
Norm and seminorm on these functional spaces are denoted by
$\NORM{\cdot}{L^{p}(\omega)}$, $\NORM{\cdot}{W^{m,p}(\omega)}$ and
$\abs{\,\cdot\,}_{W^{m,p}(\omega)}$, respectively.
For $p=2$ we prefer, as usual, the notation $H^m(\omega)$ instead
of $W^{m,2}(\omega)$, and the corresponding norm and seminorm are 
denoted by
$\NORM{\cdot}{L^{2}(\omega)}$,
$\NORM{\cdot}{H^{m}(\omega)}$ and
$\abs{\,\cdot\,}_{H^{m}(\omega)}$.
With a minor overloading of notation, we use the same symbols to
denote norms and seminorms of vector fields, e.g.,
$\NORM{\bv}{H^{m}(\omega)}$ denotes the $H^{m}$-norm of the vector
function $\bv\in\big(H^m(\omega)\big)^d$.
We denote the vector fields whose components and divergence are
in $L^2(\omega)$ by $\Hdiv(\omega)$.
We denote the space of the polynomials defined on $\omega$ of degree
$0$ and $1$ by, respectively, $\PS{0}(\omega)$ and $\PS{1}(\omega)$.
Finally, we denote the $L^2(\O)$ product between two scalar functions
$p$ and $q$ and two vector functions $\bu$ and $\bv$ by $(p,q)$ and
$(\bu,\bv)$, respectively.


\section{Mixed formulation of the diffusion problem}
\label{sec:mixed:formulation}

Let $\O\subset\Re^d$ be an open bounded polyhedral domain for $d=3$ or
a polygonal domain for $d=2$ with Lipschitz continuous boundary
$\Gamma$.
We consider the linear diffusion problem in mixed form for the scalar
unknown $p$, also dubbed the \emph{pressure}, and the vector field
$\bu$, also dubbed the \emph{pressure gradient} or, simply, the
\emph{gradient}, which reads as
\begin{align}
  \bu         &= -\nabla p \phantom{bg}\quad\textrm{in~}\O,\label{PM:a}\\
  \div(\k\bu) &= b \phantom{-\nabla pg}\quad\textrm{in~}\O,\label{PM:b}\\
  p           &= g \phantom{-\nabla pb}\quad\textrm{on~}\Gamma.\label{PM:c}
\end{align}
Hereafter, $\k(\bx)$ for $\bx\in\O$ is a possibly discontinuous, scalar
function of space; $b(\bx)$ for $\bx\in\O$ is the source term; $g(\bx)$
for $\bx\in\Gamma$ is the boundary data.
When $\k$ is discontinuous equations~\eqref{PM:a}-\eqref{PM:c} does
not have a strong solution and the solution must be understood in the
weak sense.


We assume that domain $\O$ can be split into $N_{\O}$ non-overlapping,
open and connected sub-domains $\O_i$, $i=1,\ldots,N_{\O}$, such that
$\overline{\O}=\cup_{i=1}^{N_{\O}}\overline{\O}_i$.
The diffusion coefficient may have different definitions on the
sub-domains and be discontinuous across the interfaces linking the
sub-domains.
For a proper mathematical formulation of
problem~\eqref{PM:a}-\eqref{PM:c}
we need to consider a few
assumptions on the regularity of $\k$.
Under these assumptions, it can be proved that the original continuum
problem and its variational formulation are well-posed and have a
unique and stable solution $(\bu,p)$.
We formalize these requirements as follows.

\smallskip
\begin{assumption}[\textbf{Regularity and ellipticity of the diffusion coefficient}]
  We assume that:
  
  \medskip
  \begin{description}
  \item[(K1)] $k\in L^{\infty}(\O)\cap\Pi_{i=1}^{N_{\O}}W^{1,\infty}(\O_i)$;

    \smallskip
  \item[(K2)] $\k$ is uniformly bounded from below and above almost
    everywhere in $\O$, i.e., there exists two positive constants
    $\k_*$ and $\k^*$ such that: $\kappa_*\leq\k(\bx)\leq\kappa^*$ for
    a.e.  $\bx\in\O$.
    \ENDTOPIC
  \end{description}
\end{assumption}

The normal component of flux $\k\bu$ is continuous across the 
discontinuity of $\k$ at a subdomain interface.
Hence, when the degrees of freedom are associated with the gradient
and not with the flux, a special numerical treatment of $\k$ is
required to get a convergent method.

\begin{remark}
  Our approach can be extended to the more general operator
  $\div\big(\k(p)\mathbbm{K(\bx)}\nabla p\big)$ where
  $\mathbbm{K}(\bx)$ is a diffusion tensor dependent only on the
  position $\bx$ by considering the splitting $\div\big(\k(p)\cdot)$
  and $\mathbbm{K(\bx)}\nabla p$.
  This generalization will be the topic of future works.
\end{remark}

\subsection{Mesh technicalities and diffusion coefficients}
Hereafter, we use mainly 3D notations to describe the method with a
few remarks about lower dimensions.
Let $\{\Oh\}_{h}$ be a sequence of conformal partitions of $\O$ into
non-overlapping closed polyhedral cells $c$ (polygons in two
dimensions).
Each partition $\Oh$, the \emph{mesh}, is labeled by the real
parameter $h$, whose definition is given below.
The mesh regularity assumptions on the sequence $\{\O_h\}_h$ necessary
to develop a rigorous convergence theory are presented in
Section~\ref{sec:error:estimates}.
For the moment, we only assume that mesh faces match discontinuity
interfaces of $\k$ whenever $\k$ is discontinuous in $\O$ and also
consider meshes that may contain non-convex cells and cells with
hanging nodes as those provided by local refinements, e.g., Adaptive
Mesh Refinement (AMR) techniques.
Examples of such meshes can be found
in~\cite{Eymard-Henri-Herbin-Hubert-Klofkorn-Manzini:2011:FVCA6:benchmark:proc-peer,Manzini-Russo-Sukumar:2014}.
We denote the diameter of cell $c$ by $\hc$, its boundary by $\pc$,
its volume by $|c|$, its centroid (geometric barycenter) by $\bx_c$.
The mesh size parameter is the maximum of all $\hc$.
We use the symbol $f$ for a mesh face, $|f|$ for its area (edge length
in two dimensions), $\bn_f$ for its unit normal vector whose
orientation is fixed once and for all, and $\bx_f$ for its center of
gravity (edge midpoint in two dimensions).
A mesh face can be either internal or located at the external boundary
$\Gamma$.
In the former case, we denote by $c_1$ and $c_2$ the two cells sharing
the face, so that $f\subseteq\partial c_1\cap\partial c_2$; in the
latter case, we use the notation $f\subset\Gamma$ without specifying
the unique cell to which face $f$ belongs.
We denote the space of the discontinuous functions on $\O$ whose
restriction to each cell $c$ of $\Oh$ is a constant or a linear
polynomial by, respectively, $\PS{0}(\Oh)$ and $\PS{1}(\Oh)$; for
example, $q\in\PS{0}(\Oh)$ iff $q_{|c}\in\PS{0}(c)$ for every
$c\in\Oh$.

\begin{figure}
  \centering
  \begin{overpic}[width=8cm]{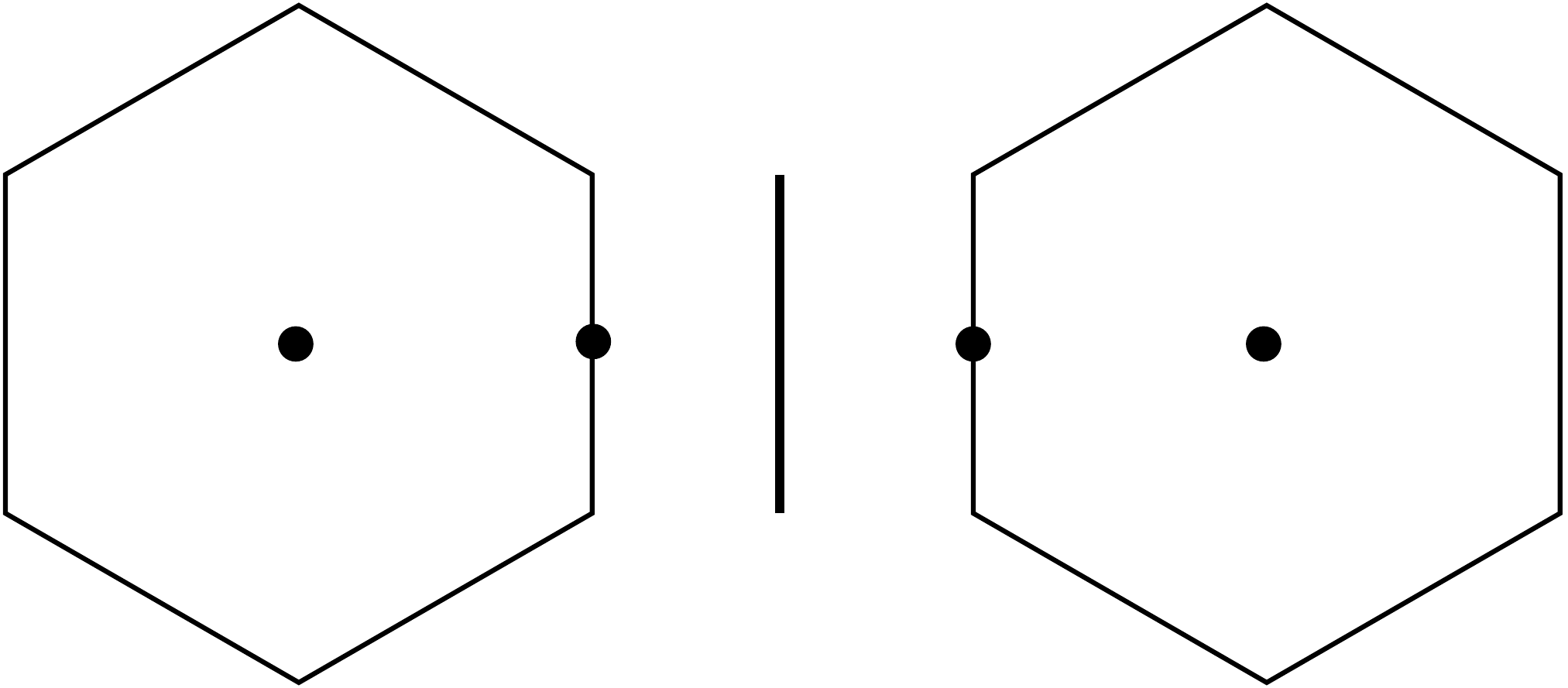}
    \put(4,25){$\k^{c_1}(\bx)\approx\k(\bx)$}
    \put(71,25){$\k^{c_2}(\bx)\approx\k$(\bx)}
    \put(30,  20){$\kcof$}
    \put(64.5,20){$\kctf$}
    \put(48,5){$f$}
    \put(43,22){$\widetilde{\k}^{c_1}_f$}
    \put(51,22){$\widetilde{\k}^{c_2}_f$}
    \put(48,35){$\widetilde{\k}_f$}
    \put(10,10){$c_1$}
    \put(85,10){$c_2$}
  \end{overpic}
  \caption{Notation for the diffusion coefficient at face
    $f\subseteq\pc_1\cap\pc_2$. For graphical convenience, the two
    cells and the common face are split. $\k^{c_i}$ is located at the
    cell-center of cell $c_i$, for $i=1$ (left cell) and $i=2$ (right
    cell); $\k_f^{c_i}$ and $\widetilde{\k}^{c_i}_f$ are associated
    with face $f$ and both refer to side $i$;
    $\widetilde{\k}_f=\widetilde{\k}^{c_1}_f=\widetilde{\k}^{c_2}_f$
    is the unique value associated with face $f$ when the two face
    coefficients coincide.}
  \label{fig:cell-notation}
\end{figure}
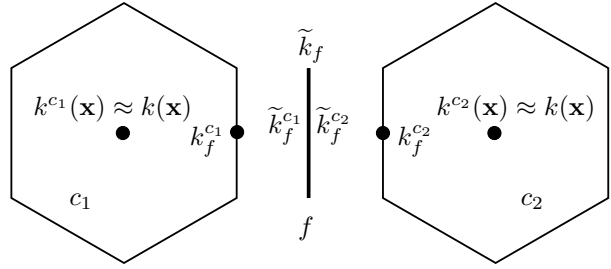

\medskip
According to Figure~\ref{fig:cell-notation}, we denote the
approximation of $\k$ associated with cell $c$ by $\kc$.
This approximation must satisfy the two following assumptions:
\begin{flalign}
  &\ASSUM{(K3)}\,\,\,\kappa_*\leq\kc(\bx)\leq\kappa^*\quad\forall\bx\in c
  \quad\textrm{and}\quad\kappa_*\leq\kcf\leq\kappa^*\quad\forall f\in\pc;   &\label{eq:kc:A}\\
  &\ASSUM{(K4)}\,\,\,\ABS{\k(\bx)-\kc(\bx)}=\kappa^*\mathcal{O}(\hc),              &\label{eq:kc:B}
\end{flalign}
where $\kappa_*$ and $\kappa^*$ are the same constants used in
\ASSUM{(K1)-(K2)}.
To this end, we may define $\kc$ as the orthogonal projection of $\k$
on either the \emph{constant} or the \emph{linear polynomials} defined
on cell $c$.
\begin{remark}
  In practice, the gradient of $\kc$ may be reconstructed from
  cell-centered values of $\k$ and may require to be limited
  to satisfy condition \ASSUM{(K3)}.
  Limiting the gradient will reduce the accuracy of the approximation
  to \ASSUM{(K4)}.
  When no limiting is used in the definition of $\kc$, it holds that
  $\NORM{\k-\kc}{L^2(c)}\leq\kappa^*\hc^{2}\ABS{\k}_{W^{1,\infty}(c)}$.
\end{remark}

Then, we introduce the \emph{one-sided face average of $\kc$ on face
  $f$}, which is given by
\begin{align}
  \kcf = \frac{1}{|f|}\Int_f\kc(x)\dS.
  \label{eq:kcf:def}
\end{align}
The values $\kcof$ and $\kctf$ provide an obvious representation of
the discontinuity of $\kc$ across face $f$.

For each internal face $f$ we introduce the \emph{face diffusion
  coefficients} $\widetilde{\k}^{c_1}_{f}$ and
$\widetilde{\k}^{c_2}_{f}$, which must satisfy the following
conditions for $i=1,2$:
\begin{flalign}
  &\ASSUM{(K5)}\,\,\,\mbox{$\widetilde{\k}^{c_i}_{f}$ only depends on $\kcof$ and $\kctf$;}                        &\label{eq:ktcf:A}\\
  &\ASSUM{(K6)}\,\,\,\kappa_*\leq\widetilde{\k}^{c_i}_{f}\leq\kappa^*;                                            &\label{eq:ktcf:B}\\
  &\ASSUM{(K7)}\,\,\,\ABS{\k^{c_i}(\bx)-\widetilde{\k}^{c_i}_{f}}\leq\kappa^*\mathcal{O}(\hc)\quad\forall\bx\in f.&\label{eq:ktcf:C}
\end{flalign}
Typically, we consider one of the two following cases:

\medskip
\begin{description}
\item[$\bullet$]
  $\widetilde{\k}_{f}:=\widetilde{\k}^{c_1}_{f}=\widetilde{\k}^{c_2}_{f}$;
  in this case $\widetilde{\k}_{f}$ is uniquely defined either as the
  arithmetic or harmonic average of $\kcof$ and $\kctf$, or by
  selecting one of the two values;

  \medskip
\item[$\bullet$] $\widetilde{\k}^{c_1}_{f}:=\kcof$ and
  $\kctf=:\widetilde{\k}^{c_2}_{f}$; in this second case two distinct
  values of the staggered diffusion coefficients are considered by
  taking the traces of $\kcf$ from the two sides of interface $f$.
\end{description}

\medskip
\begin{remark}
  In both cases, this approach preserves the symmetry and coercivity
  of the numerical formulation, thus providing a final symmetric and
  positive definite matrix operator.
\end{remark}


\smallskip
\noindent
We consider both cases in the numerical experiments of
Section~\ref{sec:numerics}.
However, since the convergence analysis only requires
conditions~\eqref{eq:ktcf:A}-\eqref{eq:ktcf:C} we do not specify the
choice of $\ktcf$ until Section~\ref{sec:numerics}.


\section{Mimetic finite difference method}
\label{sec:MFD:method}

Let $\Fh$ and $\Ph$ be the discrete spaces (formalized in
Section~\ref{subsec:discrete:spaces}) for the primary unknowns, i.e.,
pressure and pressure gradient.
Let $\buh\in\Fh$ and $\ph\in\Ph$ be the numerical approximations of
$\bu$ and $p$, respectively; $\DIVk$ the \emph{primary} mimetic
operator that approximates the \emph{combined} operator
$\div(\k\,\cdot\,)$; $\GRAD$ the \emph{derived} mimetic operator that
approximates $\nabla$; and $b^I\in\Ph$ the piecewise constant
approximation of the source term.
The definition of the discrete spaces $\Fh$ and $\Ph$, their inner
products, and the discrete divergence and gradient operators are
discussed throughout this section.

\medskip
Having introduced these quantities, the mimetic finite difference
approximation of equations \eqref{PM:a}-\eqref{PM:c} has a similar
structure and reads as:
\emph{Find $\buh\in\Fh$ and $\ph\in\Ph$ such that}
\begin{align}
  \buh       &= -\GRAD\,\ph,\label{eq:MFD-scheme:a}\\
  \DIVk \buh &= b^I.        \label{eq:MFD-scheme:b}
\end{align}
The Dirichlet boundary condition~\eqref{PM:c} are included in the
definition of the discrete differential operators $\DIVk$ and $\GRAD$
(see below).

\medskip
For $\k\in L^{\infty}(\Omega)$, $\k\bv\in\Hdiv(\O)$ and $q\in
H^1(\Omega)$, we have the integration-by-parts formula
\begin{align}
  \Int_\O\bv\cdot\k\nabla q\,\dV
  = -\Int_\O q\div\k\bv\,\dV + \Int_{\Gamma}q\bn\cdot\k\bv\,\dS.
  \label{eq:integration-by-parts}
\end{align}
Equation~\eqref{eq:integration-by-parts} implies that operator
$\div(\k\cdot)$ is in a dual relationship with the operator
$\nabla(\cdot)$.
We define the derived gradient operator by a discrete relation that
mimics~\eqref{eq:integration-by-parts}.
Let the spaces $\Fh$ and $\Ph$ be equipped with the corresponding
inner products, respectively denoted by $[\cdot,\cdot]_{\Fh}$ and
$[\cdot,\cdot]_{\Ph}$; let also $\bscalh{\cdot}{\cdot}$ be a bilinear
form that depends on boundary condition~\eqref{PM:c}.
The discrete gradient operator $\GRAD$ is derived from the primary
divergence operator $\DIVk$ according to
\begin{align}
  \scalPh{\DIVk\bvh}{\qh} = -\scalFh{\bvh}{\GRAD\qh} + \bscalh{\bvh}{g_h}
  \qquad\forall\bvh\in\Fh,\,\qh\in\Ph,
  \label{eq:duality}
\end{align}
where $g_h$ is a suitable approximation of $g=p_{|\Gamma}$.
Inclusion of boundary conditions in the definition of mimetic
operators is discussed
in~\cite{Hyman-Shashkov:1998,Shashkov-Steinberg:1996}.

\begin{remark}
  We denote the symmetric positive definite matrices representing the
  inner products in $\Fh$ and $\Ph$ by $\matM_{\Fh}$ and
  $\matM_{\Ph}$, respectively, and the matrix corresponding to the
  boundary bilinear form $\bscalh{\cdot}{\cdot}$ by
  $\matM_{\Gamma,h}$.
  Equation~\eqref{eq:duality} can be rewritten as
  \begin{align*}
    \bvh^T\matM_{\Fh}\GRAD\qh = -\bvh^T(\DIVk)^T\matM_{\Ph}\qh + \bvh^T\matM_{\Gamma,h}g_h.
  \end{align*}
  Since $\bvh$ is arbitrary, we obtain that 
  \begin{align}
    \GRAD\qh = -\matM_{\Fh}^{-1}(\DIVk)^T\matM_{\Ph}\qh + \matM_{\Fh}^{-1}\matM_{\Gamma,h}g_h.
    \label{eq:affine:grad}
  \end{align}
  Equation~\eqref{eq:affine:grad} implies that the action of $\GRAD$
  on $\qh$ is that of an affine operator, where the translation term
  depends on the Dirichlet condition $g_h$.
  Since $\qh$ is also arbitrary, when $g_h=0$, i.e., the boundary
  condition is homogeneous, equation~\eqref{eq:affine:grad} yields
  \begin{align*}
    \GRAD = -\matM_{\Fh}^{-1}(\DIVk)^T\matM_{\Ph},
  \end{align*}
  which is the matrix representation of discrete operator $\GRAD$
  in~\cite{Lipnikov-Manzini-Moulton-Shashkov:2016}.
\end{remark}


\subsection{Degrees of freedom, discrete spaces and interpolation operators}
\label{subsec:discrete:spaces}

\subsubsection{Discrete pressure space}
\label{subsubsec:discrete:pressure:space}
The members of the discrete pressure space $\cP_{h}$ consist of one
degree of freedom per cell, which represents the cell average of the
pressure.
Thus, the dimension of $\cP_{h}$ equals the number of mesh cells.
We denote the value of $\ph\in\Fh$ associated with cell $c$ by $p_c$.
Hereafter, we will conveniently identify $p_c$ with the constant
function taking this value on cell $c$ and $\ph$ with the piecewise
constant function whose restriction to cell $c$ is $p_c$.

For a given integrable scalar function $p$, we denote by $\pI\in\Ph$
the vector of degrees of freedom such that
\begin{align}
  \pI = \left\{p_c^I\right\}_{c\in\O_h},
  \qquad
  p_c^I = \Frac{1}{|c|}\Int_cp\dV.
  \label{eq:Ph:intp:I:def}
\end{align}

\subsubsection{Discrete gradient space} 
The members of the discrete gradient space $\Fth$ consist of one
degree of freedom per boundary face and two degrees of freedom per
interior face.
We denote the restriction to cell $c$ of $\buh\in\Fth$ by $\buc$ and
its component associated with face $f\in\pc$ by $\ucf$.
Hereafter, we will consider the linear subspace $\Fh$ of $\Fth$ whose
members satisfy the flux continuity constraint
\begin{equation}\label{eq:continuity}
  \widetilde{\k}^{c_1}_f\,\ucof = \widetilde{\k}^{c_2}_f\,\uctf
\end{equation}
on each internal face $f$ shared by cells $c_1$ and $c_2$.

\smallskip 
Let $\bu$ be a vector field in $\big(L^s(\O)\big)^d\cap\Hdiv(\O)$ with
$s>2$.
We define the interpolant $\buI\in\Fth$ a the vector of degrees of
freedom:
\begin{align}
  \buI = \big\{\bucI\big\}_{c \in \O_h},
  \quad
  \bucI = \big\{\big(\buI\big)^c_f\big\}_{f\in\pc},
  \quad\textrm{and}\quad
  \big(\buI\big)^c_f = \Frac{1}{|f|}\Int_f\bu_{|_c}\cdot\bn_f\dS,
  \label{eq:Fh:intp:I:def}
\end{align}
where $\bu_{|_c}$ is the restriction of $\bu$ to $c$ and
$\bu_{|_c}\cdot\bn_f$ is the one-sided limit from inside cell $c$ of
the normal component of $\bu$.


\subsection{Primary mimetic operator: the discrete divergence}
\label{subsec:DIVk}
The primary mimetic operator is the discrete divergence operator
$\DIVk\colon\Fh\to\Ph$, which is locally defined on each mesh
cell by a straightforward discretization of the divergence theorem:
\begin{align}
  \big(\DIVk \buh\big)_{|_c }
  \equiv\DIVkc\buc 
  = \frac{1}{|c|}\Sum\limits_{f\in\pc}|f|\,\scf\,\ktcf\,\ucf,
  \label{eq:DIVk:def}
\end{align}
where $\scf=\bnf\cdot\bncf$ is either $1$ or $-1$ depending on the
mutual orientation of normal $\bnf$ and the exterior normal to $\pc$
denoted by $\bncf$.
Since $\buh$ is an algebraic vector, it is convenient to think about
the discrete divergence operator as a matrix acting between the spaces
$\Fh$ and $\Ph$.
Such a matrix is full rank since $\ktcf>0$.


\subsection{Mimetic inner products and implementation}

\subsubsection{Mimetic  inner product in $\Ph$}
The mimetic inner product in space $\Ph$ is built by assembling
cell-based inner products.
Since we have only one degree of freedom per cell, this leads to a
very simple matrix representation.
The explicit formulas of the inner product in $\Ph$ are
\begin{align}
  \scalPh{\qh}{\ph} = \Sum\limits_{c \in \O_h} \scalPhc{\qh}{\ph}
  \quad
  \scalPhc{\qh}{\ph} = |c|\,p_c\,q_c.
  \label{eq:Ph:inner:product}
\end{align}
If $\qh$ and $\ph$ are the degrees of freedom of two sufficiently
regular scalar functions $q$ and $p$, i.e., $\qh=\qI$ and $\ph=\pI$,
the cell-based inner product is a second-order accurate approximation
of the $L^2(c)$ scalar product of $p$ and $q$:
\begin{align}
  \scalPhc{\qIc}{\pIc} = |c|\,\qIc\,\pIc = \Int_c p\,q\dV + |c|O(\hc^2).
  \label{eq:local:scalars}
\end{align}
Let $\matM_{\Phc}$ be the inner product matrix such that 
\begin{align}
  \scalPh{\qh}{\ph} =\qh^T\matM_{\Phc}\ph.
  \label{eq:MPh:def}
\end{align}
According to~\eqref{eq:Ph:inner:product}, $\matM_{\Phc}$ is a diagonal
matrix with values $|c|$ on the diagonal.

\subsubsection{Mimetic  inner product in $\Fh$}
The mimetic inner product in space $\Fh$ is built by assembling
cell-based inner-products to mimic the additivity of integration:
\begin{equation}\label{cell-approximation}
  \scalFh{\bvh}{\buh} 
  = \Sum\limits_{c\in\Oh}\scalFhc{\bvc}{\buc},
\end{equation}
where for every cell $c$ the local inner product
$\scalFhc{\cdot}{\cdot}$ in ${\Fh}_{|c}$ is required to satisfy the
two conditions of the following assumption.

\medskip
\begin{assumption}[\textbf{Mimetic inner product for gradients}]

  \medskip
  \begin{description}
  \item[(S1)] \textbf{spectral stability}: there exist two strictly
    positive constants $\sigma_*$ and $\sigma^*$, which are
    independent of $h$, such that for all $\buh\in\Fhc$ and for every
    cell $c$ it holds:
    \begin{align}
      \sigma_*|c|\sum_{f\in\pc}\abs{\ucf}^2\leq \scalFhc{\buc}{\buc} \leq \sigma^*|c|\sum_{f\in\pc}\abs{\ucf}^2;
      \label{eq:S1}
    \end{align}
    
  \item[(S2)] \textbf{local consistency}: for every $\buh\in\Fhc$ and
    every linear polynomial $\qone\in\PS{1}(c)$ with zero average over
    $c$ it holds:
    \begin{align}
      \scalFhc{\buc}{(\nabla\qone)^I} = 
      \sum_{f\in\pc}\scf\ucf\Int_f\kc\qone\,\dS,
      \label{eq:S2}
    \end{align}
    where $\bncf$ is the unit vector
    orthogonal to $f$ and pointing out of $c$.
  \end{description}
\end{assumption}

When $\bvh$ and $\buh$ are the degrees of freedom of two sufficiently
regular vector fields, i.e., $\bvh=\bvI$ and $\buh=\buI$, the mimetic
inner product defined by \ASSUM{(S1)-(S2)} is a local first-order
accurate approximation of the weighted $L^2(c)$ inner product of $\bv$
and $\bu$:
\begin{equation}\label{eq:cell-approximation}
  \scalFhc{\bvcI}{\bucI} = \Int_c\k\bv\cdot\bu\dV + |c|O(\hc).
\end{equation}
To prove this, we derive~\eqref{eq:S2} through a few approximation
steps.
First, we replace the vector function $\bv$ by $\bvzero$, the $L^2(c)$
orthogonal projection onto constant vectors inside $c$, which leads to
an admissible error of order $\hc$.
Second, we substitute function $\k$ with its cell-based approximation
$\kc$.
Third, we approximate function $\bu$ by a function (still denoted by
$\bu$ for simplicity of exposition) that has two special properties:
$(i)$, $\bu\cdot\bn$ is constant on each face $f$ of $\pc$, and
$(ii)$, $\div(\kc\bu)$ is constant in $c$.
The space of such functions is denoted by $\Shc$ and is sufficiently
rich to contain the constant vector functions, thus ensuring that the
approximation is convergent and (at least) first-order accurate.
Then, we show that
\begin{equation}
  \scalFhc{\bvzero^I}{\bucI} 
  = \int_c\kc\bvzero\cdot\bu\,\dV, 
  \label{eq:MFh:def}
\end{equation}
for any constant $\bvzero$ and $\bu\in\Shc$.
Since $\bvzero$ is constant on $c$, we can write $\bvzero=\nabla\qone$
where $\qone$ is a linear polynomial with zero average over $c$.
Inserting it in the right-hand side of \eqref{eq:MFh:def} and
integrating by parts, we obtain
\begin{equation}
  \label{eq:IBP_local}
  \Int_c\kc\,\nabla\qone\cdot\bu\,\dV =
  -\Int_c\div(\kc\bu)\qone\,\dV
  + \Int_{\pc}\kc\bu\cdot\bn\,\qone\,\dS.
\end{equation}
The volume integral in the right-hand side of~\eqref{eq:IBP_local} is
zero because $\div(\kc\bu)$ is assumed constant on $c$ and can be
pulled out of the integral.
By our assumptions, $\bu\cdot\bnf$ is also constant on face $f$ and
can be pulled out of the face integrals.
Since $\bu\cdot\bnf=\scf\ucf$, we obtain~\eqref{eq:S2} by defining
the inner product matrix from
\begin{equation}
  \label{eq:def:FM}
  \scalFhc{\bucI}{\bvzero} =
  ((\nabla\qone)_c^I)^T\,\matM_{\Fh,c}\,\bucI = 
  \sum\limits_{f\in\pc}\ucf\,\scf\Int_f\kc\,\qone\,\dS
  \qquad \forall\qone\in{\cal P}_1(c),\,\forall\bu\in{\cal S}(c).
\end{equation}

\begin{remark}
An important difference between this formulation and the original MFD
formulation
in~\cite{Brezzi-Lipnikov-Shashkov:2005,Brezzi-Lipnikov-Simoncini:2005}
is that in~\eqref{eq:def:FM} $\kc$ can be a linear approximation of
the diffusion coefficient $\k$.
\end{remark}

\medskip
Now, we use the linearity of the space of linear functions $\qone$ to
get an alternative representation of equation~\eqref{eq:def:FM}.
Consider the cell-based vector $\brc = \brc(\qone)$ with the following
entries:
\begin{align*}
  \brc = \{\rcf\}_{f\in\pc}, \qquad 
  \rcf = \scf\,\Int_f\kc\,\qone\,\dS.
\end{align*}

\medskip
From~\eqref{eq:def:FM}, matrix $\matM_{\Fh,c}$ is the solution of the
system of matrix equations:
\begin{equation}
  \label{matrix-equation}
  \matM_{\Fh,c} (\nabla\qone)_c^I = \brc(\qone)
  \quad \forall\qone\in\PS{1}(c).
\end{equation}
Due to linearity of these equations, it is sufficient to consider only
three linearly independent functions in 3D: $q_{1,x}=x-x_c$,
$q_{1,y}=y-y_c$, and $q_{1,z}=z-z_c$ (only $q_{1,x}$ and $q_{1,y}$ in
2D).
Let
\begin{align}
  \matNc = \big[(\nabla q_{1,x})_c^I\ (\nabla q_{1,y})_c^I\ (\nabla q_{1,z})_c^I\big],
  \quad 
  \matRc = [\brc(q_{1,x})\ \brc(q_{1,y})\ \brc(q_{1,z})]
  \label{eq:Nc:Rc:def}
\end{align}
be two column-partitioned rectangular matrices.

Matrix equation \eqref{matrix-equation} is equivalent to
\begin{equation}
  \label{MN=R}
  \matM_{\Fh,c}\matNc=\matRc.
\end{equation}

\begin{lemma}[Characterization of $\matNc^T\matRc$]
  \label{lem:NtR}
  Let matrices $\matNc$ and $\matRc$ be defined as
  in~\eqref{eq:Nc:Rc:def}, and $\kc$ be the constant or linear
  approximation of $\k$ inside cell $c$.
  Then, $\matNc^T\matRc$ is the SPD matrix given by
  \begin{align*}
    \matNc^T\matRc = \matI\,\Int_c\kc\,\dV = \kc(\bx_c)|c|\,\matI.
  \end{align*}
\end{lemma}

\medskip\noindent
\begin{proof}
  Using~\eqref{eq:def:FM}, the dot product of the first column vectors
  of matrices $\matNc$ and $\matRc$ is
  \begin{align*}
    \big((\nabla q_{1,x})_c^I\big)^T \brc(q_{1,x}) 
    = (\nabla q_{1,x})_c^I)^T \matM_{\Fh,c} (\nabla q_{1,x})_c^I
    = \Int_c\kc\nabla q_{1,x}\cdot\nabla q_{1,x}\,\dV
    = \Int_c\kc\,\dV.
  \end{align*} 
  A similar argument works for the dot products of other column
  vectors.
  The last statement of the lemma follows from exact integration of a
  constant or linear function.
\end{proof}

This lemma allows us to write matrix $\matM_{\Fh,c}$ according to the
mimetic formula~\cite{BeiraodaVeiga-Lipnikov-Manzini:2014}:
\begin{align}
  \matM_{\Fh,c} = \matRc(\matRc^T\matNc)^{-1} \matRc^T + \gamma_c\matPc,
  \quad \matPc = \matI-\matNc (\matNc^T\matNc)^{-1}\matNc^T
  \label{eq:matM:def}
\end{align}
with a positive factor $\gamma_c$ in front of the projection matrix
$\matPc$.
A recommended choice for $\gamma_c$ is the mean trace of the first
term.
A family of mimetic schemes is obtained if we replace $\gamma_c$ by an
arbitrarily symmetric positive definite matrix $\matGc$:
\begin{align*}
  \matM_{\Fh,c} = \matRc (\matRc^T\matNc)^{-1} \matRc^T + \matPc\,\matGc\,\matPc.
\end{align*}
Stability of the resulting mimetic scheme depends on spectral bounds
of matrix $\matGc$ that should be close to the value of $\gamma_c$
(see~\ASSUM{(S1)}.).

\begin{remark}
  Consider the following matrix equation
  \begin{align*}
    \matW_{\Fh,c} \matRc = \matNc.
  \end{align*}
  The solution of this equation is the inverse of matrix
  $\matM_{\Fh,c}$ for some value of $\gamma_c$ or $\matGc$.
  Only this matrix is needed in the hybridization procedure.
  The general formula for $ \matW_{\Fh,c}$ is given
  by~\cite{BeiraodaVeiga-Lipnikov-Manzini:2014}
  \begin{align*}
    \matW_{\Fh,c} 
    = \matNc(\matNc^T\matRc)^{-1}\matNc^T + \widetilde\matPc\,\widetilde\matGc\,\widetilde\matPc,
  \end{align*}
  where $\widetilde\matPc=\matI-\matRc(\matNc^T\matRc)^{-1}\matRc^T$.
\end{remark}
%


\begin{remark}[Implementation details]
  Since the diffusion coefficient $\kc$ is either a constant or a
  linear function and $\qone$ is a linear function, we can easily
  integrate $\kc\qone$ analytically or by using a sufficiently
  accurate quadrature rule (e.g., the Simpson rule on a decomposition
  of $c$ in simplexes).
  In 3D, we can also reduce the numerical integration to a 2D
  integration over the faces of $\pc$ by using the divergence theorem.
\end{remark}

\section{Stability and convergence analysis}
\label{sec:error:estimates}

In this section we prove the stability of the method (\emph{inf-sup}
condition) and the convergence of the approximation of pressure and
gradient by deriving an estimate for both errors.
To carry out the analysis of the method, we find it convenient to
reformulate~\eqref{eq:MFD-scheme:a}-\eqref{eq:MFD-scheme:b} by
using~\eqref{eq:duality} in the following pseudo-variational form:
\emph{Find $\buh\in\Fh$ and $\ph\in\Ph$ such that}
\begin{align}
  \scalFh{\buh}{\bvh} - \scalPh{\DIVk\bvh}{\ph} &= -\bscalh{\bvh}{g_h}\phantom{\scalPh{b^I}{\qh}}\qquad\forall\bvh\in\Fh,\label{eq:psvar:MFD:a}\\[0.25em]
  \scalPh{\DIVk\buh}{\qh}                       &= \,\scalPh{b^I}{\qh}\phantom{-\bscalh{g_h}{\bvh}}\qquad\forall\ph\in\Ph. \label{eq:psvar:MFD:b}
\end{align}
Formulations~\eqref{eq:MFD-scheme:a}-\eqref{eq:MFD-scheme:b}
and~\eqref{eq:psvar:MFD:a}-\eqref{eq:psvar:MFD:b} are equivalent,
except that the Dirichlet boundary conditions are now included in the
right-hand side of~\eqref{eq:psvar:MFD:a} through the term
\begin{align}
  \bscalh{\bvh}{g_h} = \sum_{f\subset\Gamma}|f|\kcf v_fg_f,
\end{align}
where $g_h=(g_f)_{f\in\Gamma}$ is the face average of $g$ on face $f$,
and $v_f$ is the value of $\bvh$ associated with $f$ (here, we omit
the superscript ``c'' as the cell is unique).
For the sake of the presentation, we consider only the case of
homogeneous Dirichlet boundary condition, i.e., $g=0$ on $\Gamma$.


The estimate of the approximation error for the gradient is carried
out in the mesh dependent norm
\begin{align*}
  \TNORM{\bvh}{\Fh}^2 
  = \Sum_{c\in\Omega_h}\TNORM{\bvc}{\Fh,c}^2
  = \Sum_{c\in\Omega_h}\scalFhc{\bvc}{\bvc},
\end{align*}
which is the norm induced by the inner product in $\Fh$.
The estimate of the approximation error for the pressure is carried
out in the mesh dependent norm
\begin{align*}
  \TNORM{\qh}{\Ph}^2 
  = \Sum_{c\in\Omega_h} \TNORM{ q_c }{ \Phc }^2
  = \Sum_{c\in\Omega_h} \scalPhc{ q_c }{ q_c }
  = \Sum_{c\in\Omega_h} |c|\,q_c^2,
\end{align*}
which is the norm induced by the inner product in $\Ph$.
Since we can identify $\qh\in\Ph$ with $\qh\in\PS{0}(\Oh)$, the
piecewise constant function defined on $\Oh$ such that
${\qh}_{|c}=q_c$, we write that
$\TNORM{\qh}{\Ph}=\NORM{\qh}{L^2(\O)}$.

As $[\cdot,\cdot]_{\Fh}$ and $[\cdot,\cdot]_{\Ph}$ are inner products,
the Cauchy-Schwarz inequalities hold:
\begin{align} 
  \label{eq:CS-global}
  \scalFh{\buh}{\bvh}\leq\TNORM{\buh}{\Fh}\TNORM{\bvh}{\Fh}
  \qquad\textrm{and}\qquad
  \scalPh{\ph}{\qh}\leq\TNORM{\ph}{\Ph}\TNORM{\qh}{\Ph}.
\end{align}
In this section we will also use the \emph{local} Cauchy-Schwarz
inequality for the gradient fields
\begin{align} 
  \label{eq:CS-cell}
  \scalFhc{\buc}{\bvc}\leq\TNORM{\buc}{\Fh,c}\TNORM{\bvc}{\Fh,c}
  \qquad\forall c\in\Oh,
\end{align}
which holds because $[\cdot,\cdot]_{\Fh,c}$ is an inner product on
${\Fh}_{|c}$.

\subsection{Mesh regularity, polynomial interpolation estimate and
  trace inequality}
\label{subsec:mesh-regularity}
The convergence analysis requires a few assumptions on the sequence of
meshes $\{\O\}_h$ that are not restrictive in practice.

\medskip
\begin{description}
\item[(MR)] There exist two positive real numbers $\Ns$ and
  $\rho_s$ such that every mesh $\{\Oh\}_{h}$ admits a conforming
  decomposition $\Th$ into shape-regular tetrahedra such that
  \begin{description}
  \item[(MR1)] every polyhedron $c$ admits a decomposition
    $\Thc$ made of less than $\Ns$ tetrahedra that includes all
    vertices of $c$;
  \item[(MR2)] each tetrahedron $T\in\Th$ is shape-regular,
    i.e., it holds that
    \begin{align}
      \rho_s\hT\leq\rT,
    \end{align}
    where $\rT$ and $\hT$ are the radius of the inscribed sphere in
    $T$ and the diameter of $T$, respectively.
  \end{description}
\end{description}
These assumptions impose some restrictions on the shape of the
admissible cell $c$ to avoid pathological situations.
Under assumption~\ASSUM{(MR)}, it is possible to prove the following
properties on the mesh, which we use in the analysis of the next
sections~\cite{BeiraodaVeiga-Lipnikov-Manzini:2014,Brenner-Scott:2002}.
Moreover, it is worth mentioning that $\Th$ is never built in the
practical implementation of the method.

\medskip
\begin{description}
\item[(M1)] The number of faces and edges of every cell $c$ is
  uniformly bounded by a constant that depends only on $\Ns$ and
  $\rho_s$.

  \smallskip
\item[(M2)] For every cell $c\in\Oh$, all the related geometric
  quantities scales in a uniform way, i.e., there exists a constant
  $\astar$ such that:
  \begin{align}
    \astar\hc^d\leq |c|\leq \hc^d\quad \forall c\in\Oh         \label{eq:M2:A}\\[0.5em]
    \astar\hc^{d-1}\leq |f|\leq \hc^{d-1}\quad \forall f\in\pc, \label{eq:M2:B}
  \end{align}
  where $d=2,3$.
  Combining~\eqref{eq:M2:A} and~\eqref{eq:M2:B} we find that 
  \begin{align}
    \frac{|c|}{|f|}\leq \frac{\hc^d}{\astar\hc^{d-1}} = \aStar\hc,
    \label{eq:M2:C}
  \end{align}
  where we set $\aStar=(\astar)^{-1}$.
  
\item[(M3)] There exists a constant $\bstar$ depending only on $\Ns$
  and $\rho_s$ such that for all $c\in\Oh$ and all $T\in\Thc$ it holds
  $\bstar\hc\leq\hT$.

\item[(M4)] \emph{(Agmon inequality)}. There exists a constant
  $\CAg>0$, which is independent of $h$, such that the following trace
  inequality, dubbed \emph{Agmon inequality}, holds true:
  \begin{align}
    \NORM{\psi}{L^2(f)}^2\leq
    \CAg\Big(
    \hc^{-1}\NORM{\psi}{L^2(c)}^2 + \hc\abs{\psi}_{H^1(c)}^2
    \Big)
    \qquad\forall f\in\pc.
    \label{eq:Agmon:inequality}
  \end{align}

\item[(M5)] \emph{(Interpolation inequalities)}. There exists a
  constant $\CIp>0$, which is independent of $h$, such that for every
  cell $c\in\Oh$ and every function $\psi\in H^2(c)$ there exists a
  constant polynomial $\psizero$ and a linear polynomial $\psi_1$
  defined on $c$ such that:
  \begin{align}
    \NORM{\psi-\psizero}{L^2(c)}                               & \leq \CIp\hc\NORM{\psi}{H^1(c)},\label{eq:Fh:intp:error:P0}\\[0.5em]
    \NORM{\psi-\psi_1}{L^2(c)} + \hc\NORM{\psi-\psi_1}{H^1(c)} & \leq \CIp\hc^2\NORM{\psi}{H^2(c)}\label{eq:Fh:intp:error:P1}.
  \end{align}
\end{description}


\subsection{Second interpolation operator and preliminary lemmas}

The interpolant defined in~\eqref{eq:Fh:intp:I:def} does not satisfy the
continuity condition~\eqref{eq:continuity} when $\k$ is discontinuous
across the mesh interface $f$.
For this reason, in the convergence analysis we need a second
interpolation operator, here denoted by $\bvII$, which is defined as
follows for gradient fields $\bv$ such that
$\k\bv\in(L^s(\O))^d\cap\Hdiv(\O)$, $s>2$, and diffusion coefficients
$\k$ satisfying assumptions~\ASSUM{(K1)-(K2)}:
\begin{align}
  \bvII = \left\{\bvcII\right\}_{c\in\O_h},
  \quad
  \bvcII = \left\{ \cf{\big(\bvII\big)} \right\}_{f \in \pc},
  \quad
  \cf{\big(\bvII\big)} = \Frac{1}{\ktcf\,|f|}\Int_f\kk_{|_c}\,(\bv\cdot\bnf)\,\dS.
  \label{eq:Fh:intp:II:def}
\end{align}
We state the properties of this second interpolation operator in the
following lemmas that are preliminary to the convergence analysis of
the next two subsections.
In all the following lemmas we assume the mesh regularity in
accordance with~\ASSUM{(MR1)-(MR2)} so that
properties~\ASSUM{(M1)-(M5)} hold.

\medskip
\begin{lemma}[Commuting property]
  \label{lemma:commuting:property}
  For every vector function $\bv$ such that $\kk\bv\in\Hdiv(\Omega)$
  it holds that
  \begin{align*}
    \DIVk\bvII = (\div(\kk\bv))^I.
  \end{align*}
\end{lemma}
\begin{proof}
  Consider a cell $c\in\Oh$.
  We use the definition of the discrete divergence operator given
  in~\eqref{eq:DIVk:def}, definitions~\eqref{eq:Fh:intp:II:def}
  and~\eqref{eq:Ph:intp:I:def} for the interpolation operators in
  $\Fh$ and $\Ph$, and we apply the Divergence Theorem to obtain:
  \begin{align*}
    \DIVkc\bvcII 
    = \frac{1}{|c|}\Int_{\pc}\kk_{|c}\bv\cdot\bn_c\,\dS
    = \frac{1}{|c|}\Int_c\div(\kk\bv)\,\dV
    = (\div(\kk\bv))^I_c,
  \end{align*}
  where $\bn_c$ is the unit vector orthogonal to $\pc$.
  The assertion of the lemma follows by collecting the relation above
  for all the cells of the mesh.
\end{proof}

\medskip
\begin{lemma}
  \label{lemma:DIVk:ortho}
  Let $\bu$ be the solution of problem \eqref{PM:a}-\eqref{PM:c},
  $\buII$ its second interpolant according
  to~\eqref{eq:Fh:intp:II:def}, and $\buh$ the discrete pressure
  gradient field solving the mimetic finite difference
  scheme~\eqref{eq:psvar:MFD:a}-\eqref{eq:psvar:MFD:b}.
  Then,
  \begin{align}
    \label{eq:DIVk:ortho}
    \DIVk\buII = \DIVk\buh.
  \end{align}
\end{lemma}
\begin{proof}
  Lemma~\ref{lemma:commuting:property} for $\bv=\bu$,
  and equations~\eqref{PM:b} and~\eqref{eq:MFD-scheme:b} imply that
  $\DIVk\buII=(\div(\kk\bu))^I=b^I=\DIVk\buh$,
  which is the assertion of the lemma.
\end{proof}

\medskip
\begin{lemma}
  \label{lemma:Fhc:norm:equivalence}
  For every vector field $\bv\in H^1(c)$ and its first and second
  interpolants $\bvcI$ and $\bvcII$ it holds that
  \begin{align}
    \TNORM{ \bvc^{I} }{\Fh,c}^2 + \TNORM{ \bvc^{\II} }{\Fh,c}^2 
    \leq\CNST{\ref{lemma:Fhc:norm:equivalence}}\,\Big(\NORM{\bv}{L^2(c)}^2+\hc^2\,\abs{\bv}_{H^1(c)}^2\Big),
    \label{eq:Fhc:norm:equivalence}
  \end{align}
  where the positive constant $\CNST{\ref{lemma:Fhc:norm:equivalence}}$ is independent of $h$.
\end{lemma}
\begin{proof}
  Inequality~\eqref{eq:Fhc:norm:equivalence} follows from the
  stability condition~\ASSUM{(S1)}, the definition of the second
  interpolant~\eqref{eq:Fh:intp:I:def}, noting that
  $\kk(\bx)\slash{\ktcf}\leq\kappa^*\slash{\kappa_*}$ for $\bx\in c$,
  applying the Agmon inequality, using~\eqref{eq:M2:C} and noting 
  that the number of faces $f\in\pc$ is uniformly bounded by $\Ns$:
  \begin{align}
    \TNORM{ \bvc^{I} }{\Fh,c}^2 + \TNORM{ \bvc^{\II} }{\Fh,c}^2 
    &\leq \sigma^*|c|\sum_{f\in\pc}\big| \cf{(\bvI) } \big|^2 
    +     \sigma^*|c|\sum_{f\in\pc}\big| \cf{(\bvII)} \big|^2
    \nonumber\\[0.5em]
    &= \sigma^*|c|\sum_{f\in\pc}\bigg(\frac{1}{|f|}\Int_f\bv\cdot\bnf\dS\bigg)^2
    +  \sigma^*|c|\sum_{f\in\pc}\bigg(\frac{1}{|f|\ktcf}\Int_f\kk\bv\cdot\bnf\dS\bigg)^2
    \nonumber\\[0.5em]
    &\leq \sigma^*\left(1+\left(\frac{\kappa^*}{\kappa_*}\right)^2\right)
    \,|c|\sum_{f\in\pc}|f|^{-1}\NORM{\bv}{L^2(f)}^2
    \nonumber\\[0.5em]
    &\leq \,\sigma^*\left(1+\left(\frac{\kappa^*}{\kappa_*}\right)^2\right)\,\CAg\,\Ns\,\aStar\,
    \Big( \NORM{\bv}{L^2(c)}^2 + \hc^2\abs{\bv}_{H^1(c)} \Big),
  \end{align}
  Finally, we set
  \begin{align*}
    \CNST{\ref{lemma:Fhc:norm:equivalence}} =
    \sigma^*\left(1+\left(\frac{\kappa^*}{\kappa_*}\right)^2\right)\CAg\,\Ns\,\aStar
  \end{align*}
  as the constant that appears in lemma's
  inequality~\eqref{eq:Fhc:norm:equivalence}.
\end{proof}

\begin{lemma}
  \label{lemma:diff:II-I}
  Consider a function $\psi\in H^2(\O)$, its piecewise polynomial
  approximation $\psione\in\PS{1}(\Oh)$ from \ASSUM{(M5)}, and denote
  by $\nabla\psione\in\big(\PS{0}(\Oh)\big)^d$, $d=2,3$, the piecewise
  constant vector such that
  $\big(\nabla{\psione}\big)_{|c}=\nabla({\psione}_{|c})$ for every
  cell $c\in\Oh$.
  Let $(\nabla\psione)^{I}$ and $(\nabla\psione)^{\II}$ be the first
  and second interpolant of $\nabla\psione$ defined
  in~\eqref{eq:Fh:intp:I:def} and~\eqref{eq:Fh:intp:II:def},
  respectively.
  Then, it holds that
  \begin{align}
    \TNORM{(\nabla\psione)^{\II}-(\nabla\psione)^{I}}{\Fh}\leq\CNST{\ref{lemma:diff:II-I}}\,h\NORM{\psi}{H^2(\O)},
  \end{align}
  where the positive constant $\CNST{\ref{lemma:diff:II-I}}$ is independent of $h$.
\end{lemma}
\begin{proof}
  Denote $\bwh=(\nabla\psione)^{\II}-(\nabla\psione)^{I}$.
  Since $\bnf\cdot\nabla\psione$ is constant, the components of $\bwh$
  are given by:
  \begin{align*}
    \cf{ (\bwh) } 
    = \cf{ \Big( (\nabla\psione)^{\II} - (\nabla\psione)^{I} \Big ) }
    = \frac{\bnf\cdot\nabla\psione}{|f|}\Int_{f}\frac{\kk-\ktcf}{\ktcf}\,\dS
  \end{align*}
  and~\eqref{eq:kc:B} and Assumption~\ASSUM{(K7)} imply that 
  \begin{align*}
    \ABS{\wcf}\leq C\frac{\kappa^*}{\kappa_*}\hc\ABS{\nabla\psione},
  \end{align*}
  where $C$ does not depend on $\hc$.
  By using the spectral stability condition~\ASSUM{(S1)}, the
  geometric inequality~\eqref{eq:M2:C}, inequality~\eqref{eq:ktcf:C},
  Agmon inequality~\eqref{eq:Agmon:inequality}, it follows that
  \begin{align}
    \TNORM{\bw_c}{\Fh,c}^2 
    &\leq \sigma^* |c|\sum_{f\in\pc}|w^c_f|^2
    \leq \aStar\,\sigma^* \hc\sum_{f\in\pc}|f| |w^c_f|^2
    \leq C\aStar\,\sigma^* \Big(\frac{\kappa^*}{\kappa_*}\Big)^2\hc^3\sum_{f\in\pc}|f|\abs{\nabla\psione}^2
    \nonumber\\[0.5em]
    &\leq C\aStar\,\sigma^* \Big(\frac{\kappa^*}{\kappa_*}\Big)^2\hc^3\sum_{f\in\pc}\NORM{\nabla\psione}{L^2(f)}^2
    \leq  C\aStar\,\sigma^* \Big(\frac{\kappa^*}{\kappa_*}\Big)^2\,\CAg\,\hc^2\ABS{\psione}_{H^1(c)}^2
    \label{eq:diff:II-I:00}
  \end{align}
  In view of~\eqref{eq:Fh:intp:error:P1} (and since
  $\hc<\textsf{diam}(\O)$) we find that
  \begin{align}
    \ABS{\psione}_{H^1(c)}
    &\leq \ABS{\psi}_{H^1(c)} + \ABS{\psione-\psi}_{H^1(c)}
    \leq \ABS{\psi}_{H^1(c)} + \CIp\hc\ABS{\psi}_{H^2(c)}
    \nonumber\\[0.5em]
    &\leq \ABS{\psi}_{H^1(c)} + \CIp\textsf{diam}(\O)\ABS{\psi}_{H^2(c)}
    \leq C^*\NORM{\psi}{H^2(c)}.
    \label{eq:diff:II-I:05}
  \end{align}
  where $C^*=\max\big(1,\CIp\textsf{diam}(\O)\big)$.
  Using this relation in the last development of~\eqref{eq:diff:II-I:00}, we find that
  \begin{align*}
    \TNORM{\bw_c}{\Fh,c}^2 
    &\leq  C\aStar\,\sigma^* \Big(\frac{\kappa^*}{\kappa_*}\Big)^2\,\CAg\big(\CIp\textsf{diam}(\O)\big)^2\,\hc^2\NORM{\psi}{H^2(c)}^2.
  \end{align*}
  The assertion of the lemma follows by adding the previous inequality
  over all the mesh cells, noting that $\hc\leq h$, and setting the
  lemma constant
  $\big(\CNST{\ref{lemma:diff:II-I}}\big)^2=C\aStar\,\sigma^*\big({\kappa^*}\slash{\kappa_*}\big)^2\,\CAg(\CIp\textsf{diam}(\O))^2$.
\end{proof}

\medskip
\begin{lemma}
  \label{lemma:scalFh:DIVk:estimate}
  Let $\bvh\in\Fh$ and $\DIVk\bvh$ its discrete divergence given
  by~\eqref{eq:DIVk:def}; let $\psi\in H^2(\O)\cap H^1_0(\O)$,
  $\psione$ its linear interpolant 
  satisfying~\eqref{eq:Fh:intp:error:P1},
  and $(\nabla\psione)^{I}$ the first interpolant of $\nabla\psione$
  defined by~\eqref{eq:Fh:intp:I:def}.
  It holds that
  \begin{align}
    \scalFh{\bvh}{(\nabla\psione)^{I}}
    = -\scalPh{\DIVk\bvh}{(\psione)^{I}} + Q_h(\bvh,\psi)
    \label{eq:scalFh:DIVk:estimate}
  \end{align}
  where the last term is bounded by the following inequality:
  \begin{align}
    \ABS{Q_h(\bvh,\psi)}\leq \CNST{\ref{lemma:scalFh:DIVk:estimate}}\,h\TNORM{\bvh}{\Fh}\,\NORM{\psi}{H^2(\O)},
    \label{eq:Qh:estimate}
  \end{align}
  and $\CNST{\ref{lemma:scalFh:DIVk:estimate}}$ is a constant independent of $h$.
\end{lemma}
\begin{proof}
  We derive~\eqref{eq:scalFh:DIVk:estimate} from the consistency
  condition \ASSUM{(S2)} with $q=\psione-\psione(\bx_c)$,
  adding and subtracting $\ktcf$,
  by noting that $\psione(\bx_c)=(\psione^I)_{|c}$,
  and using definitions~\eqref{eq:DIVk:def}
  and~\eqref{eq:Ph:inner:product} for the discrete divergence operator
  and the mimetic inner product for discrete scalar variables,
  respectively:
  \begin{align}
    &\scalFh{\bvh}{(\nabla\psione)^{I}}
    = \scalFh{\bvh}{\big(\nabla(\psione-\psione(\bx_c))\big)^{I}}
    = \sum_{c\in\Oh}\sum_{f\in\pc}\scf\vcf\int_{f}\kc\big(\psione-\psione(\bx_c)\big)\dS
    \nonumber\\[0.5em]
    &\quad= \sum_{c\in\Oh}\sum_{f\in\pc}\scf\vcf\ktcf\int_{f}\big(\psione-\psione(\bx_c)\big)\dS
    +  \sum_{c\in\Oh}\sum_{f\in\pc}\scf\vcf\int_{f}\big(\kc-\ktcf\big)\big(\psione-\psione(\bx_c)\big)\dS
    \nonumber\\[0.5em]
    &\quad= -\sum_{c\in\Oh}\psione(\bx_c)\sum_{f\in\pc}|f|\scf\vcf\ktcf + Q_h
    = -\scalPh{\DIVk\bvh}{(\psione)^{I}} + Q_h,
    \label{eq:scalFh:DIVk:estimate-00}
  \end{align}
  where 
  \begin{align}
    Q_h
    &=\sum_{c\in\Oh}\sum_{f\in\pc}\scf\vcf\ktcf\int_{f}\psione\dS
    + \sum_{c\in\Oh}\sum_{f\in\pc}\scf\vcf\int_{f}\big(\kc-\ktcf\big)\big(\psione-\psione(\bx_c)\big)\dS
    \nonumber\\[0.5em]
    &= Q_1 + Q_2.
    \label{eq:Qh:def}
  \end{align}
  
  Now, we note that $\psi$ belongs to $H^2(\O)\cap H^1_0(\O)$; hence,
  rearranging the summation on the mesh faces, using the flux
  continuity condition~\eqref{eq:continuity} and noting that
  $\scof+\sctf=0$ yield that
  \begin{align}
    \sum_{c\in\Oh}\sum_{f\in\pc}\scf\vcf\ktcf\int_f\psi\dS
    =\sum_{f\in\Oh}\big(\scof\ktcf\vcof+\sctf\ktcf\vctf\big)\int_f\psi\dS = 0.
    \label{eq:Qh:proof:00}
  \end{align}
  Therefore, we can subtract $\psi$ to the integral argument of $Q_1$,
  and using the Cauchy-Schwarz inequality, the Agmon inequality, the
  interpolation estimate~\eqref{eq:Fh:intp:error:P1} and stability
  condition \ASSUM{(S1)} we estimate this term as follows
  \begin{align}
    Q_1 
    &= \sum_{c\in\Oh}\sum_{f\in\pc}\scf\vcf\ktcf\int_{f}\big(\psi-\psione(\bx_c)\big)\dS
    \nonumber\\[0.5em]
    &\leq \kappa^*\sum_{c\in\Oh}\sum_{f\in\pc}\scf\vcf|f|^{\frac{1}{2}}\NORM{\psi-\psione}{L^2(f)}
    \nonumber\\[0.5em]
    &\leq \kappa^*\sqrt{\CAg}\sum_{c\in\Oh}\sum_{f\in\pc}\scf\vcf|f|^{\frac{1}{2}}\hc^{-\frac{1}{2}}
    \Big( \NORM{\psi-\psione}{L^2(c)} + \hc\ABS{\psi-\psione}_{H^1(c)} \Big)
    \nonumber\\[0.5em]
    &\leq \kappa^*\sqrt{\CAg}\CIp\sum_{c\in\Oh}\sum_{f\in\pc}\scf\vcf|f|^{\frac{1}{2}}\hc^{\frac{3}{2}}\ABS{\psi}_{H^2(c)}
    \nonumber\\[0.5em]
    &\leq \aStar\,\kappa^*\sqrt{\CAg}\CIp\,h\bigg(|c|\sum_{f\in\pc}\ABS{\vcf}^2\bigg)^{\frac{1}{2}}\,\bigg(\sum_{c\in\Oh}\ABS{\psi}_{H^1(c)}^2\bigg)^{\frac{1}{2}}
    \nonumber\\[0.5em]
    &\leq \aStar\,\kappa^*\sqrt{\CAg}\CIp\,h\bigg(\sum_{c\in\Oh}\TNORM{\buh}{\Fh,c}^2\bigg)^{\frac{1}{2}}\,\NORM{\psi}{H^{2}(\O)}
    \nonumber\\[0.5em]
    &\leq  \aStar\,\kappa^*\sqrt{\CAg}\CIp\,h\,\TNORM{\buh}{\Fh}\,\NORM{\psi}{H^{2}(\O)}
  \end{align}
  Term $Q_2$ can be similarly estimated by using~\eqref{eq:ktcf:C},
  the Cauchy-Schwarz inequality, the Agmon inequality, the stability
  condition \ASSUM{(S2)}, inequality~\eqref{eq:diff:II-I:05}, which
  implies that $\ABS{\psione}_{H^1(c)}\leq C\NORM{\psi}{H^2(c)}$ for
  some positive constant $C$ independent of $h$, to obtain
  \begin{align}
    Q_{2}
    &= \sum_{c\in\Oh}\sum_{f\in\pc}\scf\vcf\int_{f}\big(\kc-\ktcf\big)\,\big(\psione-\psione(\bx_c)\big)\dS
    \nonumber\\[0.5em]
    &\leq \kappa^*\,\sum_{c\in\Oh}\hc\sum_{f\in\pc}\scf\vcf|f|^{\frac{1}{2}}\NORM{\psione-\psione(\bx_c)}{L^2(f)}
    \nonumber\\[0.5em]
    &\leq \kappa^*\sqrt{\CAg}\,\sum_{c\in\Oh}\hc\sum_{f\in\pc}\scf\vcf|f|^{\frac{1}{2}}\hc^{-\frac{1}{2}}
    \bigg( \NORM{\psione-\psione(\bx_c)}{L^2(c)} + \hc\ABS{\psione}_{H^1(c)} \bigg)
    \nonumber\\[0.5em]
    &\leq \kappa^*\sqrt{\CAg}\CIp\,\sum_{c\in\Oh}\hc\sum_{f\in\pc}\scf\vcf|f|^{\frac{1}{2}}\hc^{\frac{1}{2}}\ABS{\psione}_{H^1(c)}
    \nonumber\\[0.5em]
    &\leq \kappa^*\sqrt{\CAg}\CIp\,h\,\bigg(\sum_{c\in\Oh}|c|\sum_{f\in\pc}\ABS{\vcf}^2\bigg)^{\frac{1}{2}}\,\bigg(\sum_{c\in\Oh}\ABS{\psione}_{H^1(c)}^2\bigg)^{\frac{1}{2}}
    \nonumber\\[0.5em]
    &\leq
    C\kappa^*\sqrt{\CAg}\CIp\,\,h\,\TNORM{\buh}{\Fh}\,\NORM{\psi}{H^2(\O)}
    \label{eq:J3:estimate:20}
  \end{align}
  The assertion of the lemma follows by using the above estimates of
  $Q_1$ and $Q_2$ in~\eqref{eq:Qh:def} and setting
  $\CNST{\ref{lemma:scalFh:DIVk:estimate}}=\kappa^*\sqrt{\CAg}\CIp\,(1+\aStar)$.
\end{proof}


\subsection{Well-posedness of the MFD method (inf-sup condition)}

The MFD method presented in this paper is based on a saddle-point
formulation and its well-posedness is a straightforward consequence of
the existence of a discrete \emph{inf-sup}
property~\cite{Boffi-Brezzi-Fortin:2013}.
The discrete inf-sup property is proved
in~Theorem~\ref{theorem:inf-sup} below.

\medskip
\begin{theorem}[Inf-sup condition]\label{theorem:inf-sup}
  There exists a constant $\beta_*>0$ such that for every $\qh\in\Ph$
  there exists a vector $\bvq\in\Fh$ such that:
  \begin{align}
    (i)  & \qquad\DIVk\bvq = \qh,\nonumber\\[0.5em]
    (ii) & \qquad\beta_*\TNORM{\bvq}{\Fh} \leq \TNORM{\qh}{\Ph}.
  \end{align}
  The constant $\beta_*$ is independent of $h$.
\end{theorem}

\begin{proof}
  Let $\mathbb{R}\mathbb{T}_0(\Oh)$ denote the lowest-order
  Raviart-Thomas mixed finite element space of vector-valued functions
  defined on the mesh partition $\Th$.
  From~\cite{Boffi-Brezzi-Fortin:2013} we know that there exists a
  constant $\CRTz$ independent of $h$ such that for every scalar
  function $\qh\in L^2(\O)$ there exists a vector function
  $\bvRT\in\mathbb{R}\mathbb{T}_0(\Oh)$ that satisfies
  \begin{align}
    \div\,\big(\bvRT\big) &= \qh\quad\textrm{in~}\O\label{eq:infsup:1}\\[0.5em]
    \NORM{\bvRT}{L^2(\O)} + \NORM{\div\,\bvRT}{L^2(\O)} &\leq \CRTz\,\NORM{\qh}{L^2(\O)}.\label{eq:infsup:2}
  \end{align}
  Consider the discrete field $\bvq=\big(\k^{-1}\bvRT\big)^{\II}\in\Fh$.
  Assertion $(i)$ follows immediately since on each cell $c\in\Oh$
  Lemma~\ref{lemma:commuting:property} and
  equation~\eqref{eq:infsup:1} imply that:
  \begin{align}
    \DIVkc\bvq = \DIVkc\,\big((\k^{-1}\bvRT)^{\II}\big)_c 
    = \big(\div(\bvRT)\big)^I = (\qh)^I_{|c} = (\qh)_{|c}.
  \end{align}
  To prove assertion $(ii)$, we use
  Lemma~\ref{lemma:Fhc:norm:equivalence}, the local inverse inequality
  $\hc\NORM{\bvRT}{H^1(c)}\leq C\NORM{\bvRT}{L^2(c)}$ for some
  positive constant $C$ independent of $h$, and
  inequality~\eqref{eq:infsup:2} to obtain:
  \begin{align}
    \TNORM{\bvq}{\Fh}^2 
    &=      \sum_{c\in\Oh} \TNORM{\big(\k^{-1}\bvRT\big)^{\II}}{\Fh,c}^2
    \leq \CNST{\k}\CNST{\ref{lemma:Fhc:norm:equivalence}}\,\sum_{c\in\Oh} \Big( \NORM{\bvRT}{L^2(c)}^2 + \hc^2\NORM{\bvRT}{H^1(c)}^2 \Big)
    \nonumber\\[0.5em]
    &\leq \CNST{\k}\CNST{\ref{lemma:Fhc:norm:equivalence}}\,\sum_{c\in\Oh} \NORM{\bvRT}{L^2(c)}^2
    \leq \CNST{\k}\CNST{\ref{lemma:Fhc:norm:equivalence}}\CRTz\,\NORM{\qh}{L^2(\O)}^2,
  \end{align}
  where
  $\CNST{\k}=\max(\kappa_*^{-2},\max_{c\in\Oh}|\k^{-1}|_{W^{1,\infty}(c)}^2)$.
  The second assertion of the lemma follows from the identification of
  $\Ph$ and $\PS{0}(\Oh)$ which implies that
  $\NORM{\qh}{L^2(\O)}=\TNORM{\qh}{\Ph}$, and setting
  $(\beta_*)^{-2}=\CNST{\k}\CNST{\ref{lemma:Fhc:norm:equivalence}}\CRTz$.
\end{proof}

\subsection{Convergence estimate for the gradient}
The main result of this section is the following theorem.
\begin{theorem}
  \label{thm:error-gradient}
  Let $(p,\,\bu)$ be the solution of problem \eqref{PM:a}-\eqref{PM:c}
  under Assumption \ASSUM{(K1)-(K2)} with $g=0$, $p\in H^2(\O)\cap
  H^1_0(\O)$ and $\bu=\nabla p\in H^1(\Omega)$.
  Let $(\ph,\,\buh)\in\Ph\times\Fh$ be the solution of the mimetic
  problem \eqref{eq:psvar:MFD:a}-\eqref{eq:psvar:MFD:b} under
  Assumptions~\ASSUM{(K1)-(K7)},
  \ASSUM{(S1)-(S2)}, \ASSUM{(MR1)-(MR2)}.
  Then, it holds that
  \begin{align}
    \TNORM{ \bu^{\II} - \buh}{\Fh}\leq{\cal C}\,h
    \NORM{\bu}{H^1(\O)}, 
    \label{eq:thm:error-gradient}
  \end{align}
  where the positive constant ${\cal C}$ is independent of $h$.
\end{theorem}
\begin{proof}
  Let $\bepsh = \bu^{\II} - \buh$.
  Let $\pone$ be the piecewise linear interpolant of $p$ in
  $\PS{1}(\Oh)$ that is defined in each cell $c$ according to
  \ASSUM{(M5)}, and consider the piecewise constant vector
  $\buzero\in\PS{0}(\Oh)$ that is locally defined by
  ${\buzero}_{|c}=\nabla({\pone}_{|c})$ for each $c\in\Oh$.
  Adding and subtracting $\buzero^{\II}$ yields:
  \begin{align}
    \TNORM{ \bepsh }{\Fh}^2 = 
    \sum\limits_{c\in\Omega_h} \Big(
    \scalFhc{(\bu - \buzero)^{\II}_c}{\bepsc} + \scalFhc{(\buzero)^{\II}_c}{\bepsc}
    - \scalFhc{\buc}{\bepsc} \Big)
    = T_1 + T_2 + T_3.
    \label{eq:gradient:error}
  \end{align}
  We will estimate the three terms $T_1$, $T_2$, $T_3$ separately.

  \bigskip
  \textbf{Estimate of $T_1$.}
  Term $T_{1}$ is bounded by applying the Cauchy-Schwarz
  inequality~\eqref{eq:CS-cell} and the result of
  Lemma~\ref{lemma:Fhc:norm:equivalence} to each cell-wise component
  of $T_1$:
  \begin{align*}
    \abs{T_1}
    &\leq\sum_{c\in\Oh}\abs{\scalFhc{(\bu - \buzero)^{\II}_c}{\bepsc}}
    \leq \sum_{c\in\Oh}\TNORM{(\bu - \buzero)^{\II}_c}{\Fh,c}\,\TNORM{\bepsc}{\Fh,c}
    \nonumber\\[0.5em]
    &\leq 
    \CNST{\ref{lemma:Fhc:norm:equivalence}}\,
    \sum_{c\in\Oh}\Big(\NORM{\bu - \buzero}{L^2(c)}^2+\hc^2\,\abs{\bu}_{H^1(c)}^2\Big)^{\frac{1}{2}}\,\TNORM{\bepsc}{\Fh,c},
  \end{align*}
  and then applying the polynomial interpolation
  estimate~\eqref{eq:Fh:intp:error:P0} to obtain
  \begin{align*}
    \abs{T_1}\leq \,h\,\abs{\bu}_{H_1(\Omega)}\,\TNORM{\bepsh}{\Fh},
  \end{align*}
  where $C=\CNST{\ref{lemma:Fhc:norm:equivalence}}\CIp$.

  \bigskip
  \textbf{Estimate of $T_2$.}
  To estimate term $T_2$, we introduce the discrete field
  $\bwh=(\bw_c)_{c\in\Oh}\in\Fth$ such that
  $$\bwh=\buzero^{\II}-\buzero^{I}=(\nabla\pone)^{\II}-(\nabla\pone)^{I}.$$
  Therefore, we have
  \begin{align}
    T_2 
    = \sum_{c\in\Oh}\scalFhc{\bepsc}{(\nabla\pone)^{\II}_c} 
    = \sum_{c\in\Oh}\Big(\scalFhc{\bepsc}{(\nabla\pone)^I_c} + \scalFhc{\bepsc}{\bw_c} \Big)
    = T_{21} + T_{22}.
    \label{eq:grad:convg:00}
  \end{align}
  We bound term $T_{21}$ by using
  Lemma~\ref{lemma:scalFh:DIVk:estimate} with $\psi=p$, $\bvh=\bepsh$,
  and noting that $\DIVk\bepsh=0$ from Lemma~\ref{lemma:DIVk:ortho}:
  \begin{align}
    \ABS{T_{21}} = \ABS{ \scalFh{\bepsh}{(\nabla\pone)^I} } = \ABS{Q_h(\bepsh,p)} \leq C\,h\TNORM{\bepsh}{\Fh}\,\NORM{p}{H^2(\O)},
  \end{align}
  To estimate term $T_{22}$, we apply the Cauchy-Schwarz inequality
  and Lemma~\ref{lemma:diff:II-I} with $\psi=p$:
  \begin{align*}
    \ABS{T_{22}}
    \leq \TNORM{\bepsh}{\Fh}\,\TNORM{\bwh}{\Fh}
    \leq \CNST{ \ref{lemma:diff:II-I} }\,h\,\TNORM{\bepsh}{\Fh}\,\NORM{p}{H^2(\O)}.
  \end{align*}
  
  \bigskip
  \textbf{Estimate of $T_3$.}
  Finally, term $T_3$ is zero because Lemma~\ref{lemma:DIVk:ortho}
  implies that $\DIVk\bepsh=\DIVk(\bu^{\II}-\buh)=0$ and from
  equation~\eqref{eq:psvar:MFD:a} with $\bvh=\bepsh$ (recall that
  $g_h=0$) we have that
  \begin{align*}
    T_3 
    = \sum\limits_{c\in\Omega_h}\scalFhc{\buc}{\bepsc}
    = \scalFh{\buh}{\bepsh} 
    = \scalPh{\DIVk\bepsh}{\ph} = 0.
  \end{align*}
  Collecting the estimates for $T_1$ and $T_2$
  in~\eqref{eq:gradient:error} proves the assertion of the theorem.
\end{proof}

An immediate consequence of Theorem~\ref{thm:error-gradient} is the
convergence result for the flux approximation, which we state in the
following corollary.
\begin{corollary}
  \label{coro:error-flux} 
  Let $\matKc$ be the cell-based diagonal matrix formed by coefficients
  $\ktcf$, $f \in \pc$. 
  Under the same assumptions of Theorem~\ref{thm:error-gradient}, it
  holds that
  \begin{align}
    \left(\sum_{c\in\Oh}|c|\,\| \matKc(\buIc - \buc) \|^2\right)^{\frac{1}{2}}\leq C\,h
    \NORM{\bu}{H^1(\O)}, 
    \label{eq:coro:error-flux}
  \end{align}
  where $\|\cdot\|$ is the Euclidean norm for vectors, and the
  positive constant $C$ is independent of $h$ but may depend on
  the ellipticity constant $\kappa^*$ introduced in Assumption
  \ASSUM{(K2)}.
\end{corollary}
\begin{proof}
  The spectral equivalence stated by~\eqref{eq:S1} and
  Assumption~\ASSUM{(K2)} implies the equivalence of the left-hand
  side of~\eqref{eq:coro:error-flux} and $\TNORM{ \bu^{\II} - \buh}{\Fh}$,
  which is the left-hand of~\eqref{eq:thm:error-gradient}.
  This norm equivalence implies the assertion of the corollary.
\end{proof}

\subsection{Convergence estimate for the pressure}

In this section we prove the convergence of the pressure approximation
and derive an estimate for the approximation error.
The result of this section is stated in the following theorem.

\medskip
\begin{theorem}
  \label{theorem:pressure:convergence}
  Let $(p,\,\bu)$ be the solution of continuum problem
  \eqref{PM:a}-\eqref{PM:c} in the $H^2$-regular domain $\O$ under
  Assumption \ASSUM{(K1)-(K2)} with $g=0$, $b\in H^1(\O)$, $p\in
  H^2(\O)\cap H^1_0(\O)$ and $\bu=-\nabla p\in H^1(\Omega)$.
  Let $(\ph,\,\buh)\in\Ph\times\Fh$ be the solution of the mimetic
  problem \eqref{eq:psvar:MFD:a}-\eqref{eq:psvar:MFD:b} under
  Assumptions~\ASSUM{(K1)-(K7)}, \ASSUM{(S1)-(S2)},
  \ASSUM{(MR1)-(MR2)}.
  Then,
  \begin{align}
    \TNORM{\ph-\pI}{\Ph}\leq Ch\big(\NORM{p}{H^2(\O)} + \NORM{b}{H^1(\O)}\big),
    \label{eq:thm:pressure:convergence}
  \end{align}
  where the positive constant $C$ is independent of $h$.
\end{theorem}
\begin{proof}
\noindent
  Let $\psi$ be the solution of the auxiliary elliptic problem:
  \begin{align}
    -\div(\k\nabla\psi) &= \pI-\ph \phantom{0}\quad\textrm{in~}\O,    \label{eq:H2:reg-a}\\
    \psi                &= 0 \phantom{\pI-\ph}\quad\textrm{on~}\Gamma.\label{eq:H2:reg-b}
  \end{align} 
  We assume the \emph{$H^2$-regularity of solution $\psi$}: there
  exists a constant $C^*_{\O}>0$, which is independent of $h$ but may
  depend on the shape of domain $\O$, such that
  \begin{align}
    \NORM{\psi}{H^2(\O)}\leq C^*_{\O}\NORM{\pI-\ph}{L^2(\O)} = C^*_{\O}\TNORM{\pI-\ph}{\Ph}.
    \label{eq:H2:regularity}
  \end{align}
  We take $\bvpsi=(\nabla\psi)^{\II}$.
  From a straightforward calculation using the commutation property
  from Lemma~\ref{lemma:commuting:property} and the fact that
  $\div(\k\nabla\psi)$ is piecewise constant on $\Oh$ it follows that:
  \begin{align}
    \DIVk\bvpsi 
    = \DIVk(\nabla\psi)^{\II} = \big(\div(\k\nabla\psi)\big)^{I} 
    = \div(\k\nabla\psi) = \ph-\pI.
    \label{eq:vpsi:def}
  \end{align}
  In view of~\eqref{eq:vpsi:def}, we have
  \begin{align*}
    \begin{array}{rll}
      &\TNORM{\ph-\pI}{\Ph}^2 = \scalPh{\DIVk\bvpsi}{\ph-\pI}                                                       & \mbox{\big[ use~\eqref{eq:psvar:MFD:a} with $g_h=0$\big]}\\[1.0em]
      {}                      &\qquad= \scalFh{\buh}{\bvpsi} - \scalPh{\DIVk\bvpsi}{\pI}                            & \mbox{\big[ use~\eqref{eq:vpsi:def},~\eqref{eq:Ph:intp:I:def}, and~\eqref{eq:Ph:inner:product}\big]}\\[1.0em]
      {}                      &\qquad= \scalFh{\buh}{\bvpsi} - \sum_{c\in\Oh}\Int_c\,p\,\div(\k\nabla\psi)\,\dV      & \mbox{\big[ integrate by parts and use~\eqref{eq:H2:reg-b}\big]}\\[1.0em]
      {}                      &\qquad= \scalFh{\buh}{\bvpsi} + \sum_{c\in\Oh}\Int_c\,\k\nabla p\cdot\nabla\psi\,\dV  & \mbox{\big[ integrate by parts; use~\eqref{PM:a}-\eqref{PM:b},~\eqref{eq:H2:reg-b}\big]}\\[1.0em]
      {}                      &\qquad= \scalFh{\buh}{\bvpsi} + \Int_{\O}b\psi\dV.
    \end{array}
  \end{align*}
  Let $\psione$ be the piecewise linear interpolant that
  satisfies~\eqref{eq:Fh:intp:error:P1} on every cell $c$.
  We substitute $\bvpsi=(\nabla\psi)^{\II}$, add and subtract
  $(\nabla\psione)^{\II}$ and $(\nabla\psione)^{I}$ to obtain:
  \begin{align*}
    \scalFh{\buh}{\bvpsi} 
    &= \scalFh{\buh}{(\nabla\psi)^{\II}} 
    = \scalFh{\buh}{(\nabla\psione)^{\II}} + \scalFh{\buh}{(\nabla(\psi-\psione))^{\II}}\\[0.5em]
    &= \scalFh{\buh}{(\nabla\psione)^{I}} 
    + \scalFh{\buh}{(\nabla(\psi)^{\II}-(\nabla\psione))^{I}}
    + \scalFh{\buh}{(\nabla(\psi-\psione))^{\II}}.
  \end{align*}
  From this development it follows that
  \begin{align}
    \TNORM{\ph-\pI}{\Ph}^2 = J_1 + J_2 + J_3,
  \end{align}
  where
  \begin{align}
    J_1 &= \scalFh{\buh}{(\nabla\psione)^{I}} + \Int_{\O}b\psi\dV ,\label{eq:J1:def}\\[0.5em]
    J_2 &= \scalFh{\buh}{(\nabla\psione)^{\II} - (\nabla\psione)^{I}},\label{eq:J2:def}\\[0.5em]
    J_3 &= \scalFh{\buh}{(\nabla(\psi-\psione))^{\II}}. \label{eq:J3:def}
  \end{align}
  The proof of the theorem continues with the estimate of these three
  terms.

  \bigskip
  \textbf{Estimate of $J_1$.}
  We first transform term $J_{1}$ as follows by using
  Lemma~\ref{lemma:scalFh:DIVk:estimate} and
  equation~\eqref{eq:MFD-scheme:b}
  \begin{align}
    J_1
    = -\scalPh{\DIVk\buh}{(\psione)^{I}} + \Int_{\O}b\psi\dV + Q_h(\buh,\psi)
    = \sum_{c\in\Oh}\int_{c}\big(b\psi - b^I_{c}(\psione)^{I}_{c}\big)\dV + Q_h(\buh,\psi)
    \label{eq:J1:estimate:00}
  \end{align}
  Then, we use the interpolation estimates for the integral term and
  use the estimate of $Q_h$ provided by
  Lemma~\ref{lemma:scalFh:DIVk:estimate} and we obtain:
  \begin{align}
    \ABS{J_1} \leq C_{1}\,\,h\big( \ABS{b}_{H^1(\O)} + \TNORM{\buh}{\Fh} \big)\,\NORM{\psi}{H^2(\O)},
    \label{eq:J1:estimate}
  \end{align}
  where the final constant $C_{1}$ depends on $\CIp$ and
  $\CNST{\ref{lemma:scalFh:DIVk:estimate}}$ but is independent of $h$.

  \bigskip
  \textbf{Estimate of $J_2$.}  
  We use the Cauchy-Schwarz inequality and Lemma~\ref{lemma:diff:II-I}
  and we have that
  \begin{align}
    \abs{J_2}
    \leq \TNORM{\buh}{\Fh}\,\TNORM{(\nabla\psione)^{\II}-(\nabla\psione)^{I}}{\Fh}
    \leq C^*_2\,h\,\TNORM{\buh}{\Fh}\,\NORM{\psi}{H^2(\O)},
    \label{eq:J2:estimate}
  \end{align}
  where we only need to set $C^*_2=\CNST{\ref{lemma:diff:II-I}}$, which is 
  independent of $h$.

  \bigskip
  \textbf{Estimate of $J_3$.}
  Applying the Cauchy-Schwarz inequality to~\eqref{eq:J3:def} readily
  gives:
  \begin{align}
    \abs{J_3} \leq \TNORM{\buh}{\Fh}\,\TNORM{(\nabla(\psi-\psione))^{\II}}{\Fh}.
    \label{eq:J3:estimate:00}
  \end{align}
  We use \ASSUM{(S1)} (spectral stability) to estimate the second term
  in~\eqref{eq:J3:estimate:00}:
  \begin{align}
    \TNORM{ (\nabla(\psi-\psione))^{\II} }{\Fh}^2
    \leq \sigma^*\sum_{c\in\Oh}|c|\sum_{f\in\pc}\abs{ \cf{\Big(\big(\nabla(\psi-\psione)\big)^{\II}\Big)} }^2.
    \label{eq:J3:estimate:05}
  \end{align}
  From~\eqref{eq:Fh:intp:II:def}, \ASSUM{(K2)} and~\eqref{eq:ktcf:A}
  it follows that
  \begin{align*}
    \abs{ \cf{\Big(\big( \nabla(\psi-\psione)\big)^{\II}\Big)} }^2 
    &= \abs{ \frac{1}{|f|}\Int_f\bnf\cdot\frac{\kk}{\ktcf}\nabla(\psi-\psione)\dS }^2
    \leq 
    \bigg(\frac{\kappa^*}{\kappa_*}\bigg)^2
    \frac{1}{|f|}\Int_f\abs{\nabla(\psi-\psione)}^2\dS.
     \nonumber\\[0.5em]
     &= \bigg(\frac{\kappa^*}{\kappa_*}\bigg)^2
     \frac{1}{|f|}\NORM{\nabla(\psi-\psione)}{L^2(f)}^2.
  \end{align*}
  Then, we use Agmon inequality~\eqref{eq:Agmon:inequality} and the
  polynomial interpolation estimate~\eqref{eq:Fh:intp:error:P1} to obtain:
  \begin{align}
    \abs{ \cf{\Big(\big( \nabla(\psi-\psione)\big)^{\II}\Big)} }^2
    &\leq
    \bigg(\frac{\kappa^*}{\kappa_*}\bigg)^2\frac{1}{|f|}
    \CAg\bigg(\hc^{-1}\NORM{\nabla(\psi-\psione)}{L^2(c)}^2 + \hc\abs{\nabla(\psi-\psione)}_{H^1(c)}^2\bigg)
    \nonumber\\[0.5em]
    &\leq
    \bigg(\frac{\kappa^*}{\kappa_*}\bigg)^2\CAg(\CIp)^2\,
    \frac{\hc}{|f|}\NORM{\psi}{H^2(c)}.
    \label{eq:J3:estimate:10}
  \end{align}
  We substitute~\eqref{eq:J3:estimate:10} in the right-hand side
  of~\eqref{eq:J3:estimate:05}, and use the resulting inequality 
  in~\eqref{eq:J3:estimate:00}.
  In view of~\eqref{eq:M2:C}, we have the final bound of $J_3$, which
  reads as
  \begin{align}
    \abs{J_3} \leq
    \sigma^*\bigg(\frac{\kappa^*}{\kappa_*}\bigg)^2\CAg(\CIp)^2\,
    \TNORM{\buh}{\Fh}\,\bigg(\sum_{c\in\Oh}\sum_{f\in\pc}\frac{|c|\hc}{|f|}\NORM{\psi}{H^2(c)}^2\bigg)^{\frac{1}{2}}
    \leq C_3^*\,h\,\TNORM{\buh}{\Fh}\,\NORM{\psi}{H^2(\O)},
    \label{eq:J3:estimate}
  \end{align}
  where we set
  $C_3^*=\sqrt{\aStar\Ns}\sigma^*(\kappa^*\slash{\kappa_*})^2\CAg(\CIp)^2$,
  and note that this constant is independent of $h$.

  \medskip 
  The assertion of the theorem follows from
  estimates~\eqref{eq:J1:estimate}, \eqref{eq:J2:estimate},
  \eqref{eq:J3:estimate}, the $H^2$-regularity
  bound~\eqref{eq:H2:regularity}, and using
  Theorem~\ref{theorem:inf-sup} and
  Lemma~\ref{lemma:Fhc:norm:equivalence}.
\end{proof}


\section{Numerical experiments}
\label{sec:numerics}

Consider the scalar field $p$ and the diffusion tensor $\k$ given by
\begin{align}
  p(x,y) = 
  \begin{cases}
    a_1 x^2+y^2                       & x<0.5,\\[0.5em]
    a_2 x^2+y^2 +\frac{1}{4}(a_1-a_2) & x>0.5,
  \end{cases}
  \quad
  \k(x,y) =
  \begin{cases}
    b_1( 1 + x\sin(y) )     & x<0.5,\\[0.5em]
    b_2( 1 + 2x^2\sin(y) )  & x>0.5.
  \end{cases}
\end{align}
where $a_i,b_i$ are real constant numbers such that $a_ib_i=1$,
$i=1,2$, and $\bu=-\nabla p$.
We consider two test cases with, respectively, a continuous and a
discontinuous function $\k$.
In the first test case, we set $b_1 = b_2 = 1$.
In the second test case, we set $b_1 = 1$ and $b_2 = 20$, so that the
normal component of $\k\nabla p$ is continuous across the interface
boundary $x = 0.5$ while the tangential component is discontinuous.
We solve problem~\eqref{PM:a}-\eqref{PM:c} on $\O=\O_1\cup\O_2$, where
$\O_1=[0,0.5]\times[0,1]$ and $\O_2=[0.5,1]\times[0,1]$, using two
different realizations of the new MFD method, hereafter labeled as
``\textsf{Trace}'' and ``\textsf{Upwind}''.
In \SchemeI, the face coefficient $\kcf$ is the trace of
$\kc$.
In \SchemeII, the face coefficient $\kcf$ for all faces where
$\k$ is continuous is selected between $\kcof$ and $\kctf$ (we recall
that $f\subseteq\pc_1\cap\pc_2$) by taking the one from the cell whose
centroid has the bigger x-coordinate.
At faces where $\k$ is discontinuous, $\kcf$ is the trace of $\kc$ as
for \SchemeI.
The selection strategy of \SchemeII{} simulates the upwinding
between $\kcof$ and $\kctf$.
Note that a truly upwind strategy must follow in some sense the
``\emph{flow of information on the grid}'' and requires some
\emph{knowledge of the approximate solution}.
For this reason, upwinding is easily implementable in time-dependent
problems or non-linear problems where the solution at the previous
timestep or at the previous iteration is available.
In stationary linear problems, upwinding cannot be implemented without
introducing a non-linearity in the numerical formulation.
To avoid this collateral effect, we select one of the two face
coefficients according to a simple geometric criterion, which is
sufficient for our purpose.

\begin{figure}
  \centering
  \includegraphics[width=4cm]{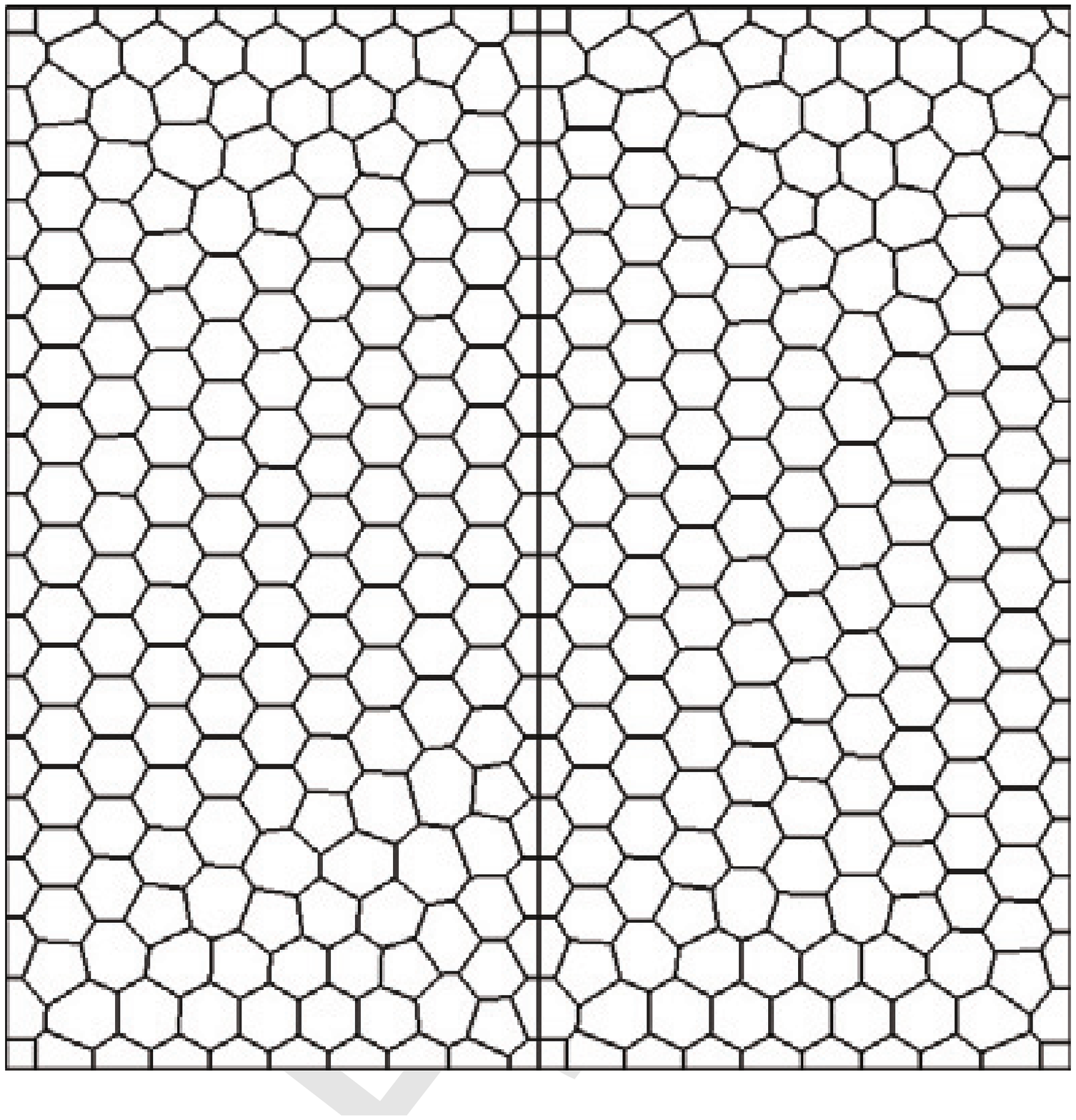}
  \vspace{-2\baselineskip}
  \caption{First mesh of the polygonal mesh sequence.}
  \label{fig:polygonal:mesh}
\end{figure}

Numerical experiments are carried out on a sequence of polygonal
meshes partitioning the two subdomains $\O_1$ and $\O_2$, see
Figure~\ref{fig:polygonal:mesh}.
To build each polygonal mesh we first generate two matching Delaunay
meshes in the left and right parts of $\O$ and, then, we build a
constraint Voronoi tessellation in each subdomain.

The relative errors for pressure and flux reads as:
\begin{align*}
  err(p) = \Frac{|||\pI - \ph |||_{\Ph}}{||| \pI |||_{\Ph}};
  \qquad
  err(\k\bu) = 
  \Frac
  {\left(\sum_{c\in\Oh}|c|\| \matKc (\buIc - \buc) \|^2\right)^{\frac{1}{2}}}
  {\left(\sum_{c\in\Oh}|c|\| \matKc \buI \|^2\right)^{\frac{1}{2}}}.
\end{align*} 
where $\pI$ is the interpolation of the exact solution $p$ defined
in~\eqref{eq:Ph:intp:I:def}, $\buI$ is the first interpolation of the
exact solution gradient $\bu=-\nabla p$ defined
in~\eqref{eq:Fh:intp:I:def}; $\matKc$ is the cell-based diagonal
matrix formed by coefficients $\ktcf$ introduced in
Corollary~\ref{coro:error-flux}.
%
For quasi-uniform meshes considered in the numerical experiments, the
Euclidean norm leads to the same conclusions as any reasonable
mesh-dependent $L^2$ norm.

We report the approximation errors when $\k$ is continuous in
Table~\ref{fig:continuous:K} and when $\k$ is discontinuous in
Table~\ref{fig:discontinuous:K}.
When $\kc$ is piecewise constant, the convergence rate of \SchemeII{}
agrees with the theory since the approximation errors of pressure and
flux scale down linearly as expected from
estimate~\eqref{eq:thm:pressure:convergence} of
Theorem~\ref{theorem:pressure:convergence} and estimate
\eqref{eq:coro:error-flux} in Corollary~\ref{coro:error-flux}.
In the other cases, a superconvergence effect is visible as the
pressure approximation rate is close to $h^{2}$ and the velocity
approximation rate is close to $h^{\frac{3}{2}}$.
Accordingly, \SchemeI{} is more accurate than
\SchemeII{} when $\kc$ is piecewise constant, while the
accuracy of these schemes is almost the same in the other cases.
The superconvergence effect will be investigated in a future work.
Finally, the behavior of \SchemeI{} and \SchemeII{} is
essentially the same regardless of $\k$ being continuous or
discontinuous.

\begin{table}
  \centering
  \caption{Relative approximation errors and convergence rates for $p$ 
    and $\bu$ in the case of a continuous diffusion coefficient $\k$ 
    using polygonal meshes as the one shown in Figure~\ref{fig:polygonal:mesh}.}
  \label{fig:continuous:K}
  \begin{tabular}{c||cc|cc||cc|cc}
    &\multicolumn{4}{c||}{$\kc\in\PS{0}(\Oh)$}
    &\multicolumn{4}{c}  {$\kc\in\PS{1}(\Oh)$}\\
    \hline
    \textsf{cells}
    &\multicolumn{2}{c|}{\SchemeI}&\multicolumn{2}{c||}{\SchemeII}
    &\multicolumn{2}{c|}{\SchemeI}&\multicolumn{2}{c}  {\SchemeII}\\
    \hline
    {} 
    &\textsf{err(p)} & \textsf{err(\k\bu)} & \textsf{err(p)} & \textsf{err(\k\bu)}
    &\textsf{err(p)} & \textsf{err(\k\bu)} & \textsf{err(p)} & \textsf{err(\k\bu)}\\
    \hline
    412     & 3.220e-3 & 7.088e-3 & 7.840e-3 & 3.877e-2 & 2.621e-3 & 3.029e-3 & 2.629e-3 & 3.050e-3\\
    1591    & 7.913e-4 & 2.251e-3 & 4.440e-3 & 1.967e-2 & 6.442e-3 & 9.619e-4 & 6.450e-4 & 9.656e-4\\
    6433    & 1.904e-4 & 8.406e-4 & 2.627e-3 & 9.844e-3 & 1.544e-3 & 4.407e-4 & 1.544e-4 & 4.413e-4\\
    25698   & 4.716e-5 & 2.432e-4 & 1.360e-3 & 4.818e-3 & 3.817e-5 & 1.314e-4 & 3.819e-5 & 1.314e-4\\
    102772  & 1.167e-5 & 1.123e-4 & 6.320e-4 & 2.531e-3 & 9.513e-6 & 5.687e-5 & 9.515e-6 & 5.688e-5\\
    \hline
    \textsf{rate} & 2.03 & 1.52 & 0.90 & 0.99 & 2.37 & 1.44 & 2.04 & 1.44
  \end{tabular}
\end{table}

\begin{table}
  \centering
  \caption{Relative approximation errors and convergence rates for $p$ 
    and $\bu$ in the case of a discontinuous diffusion coefficient $\k$ 
    using polygonal meshes as the one shown in Figure~\ref{fig:polygonal:mesh}.}
  \label{fig:discontinuous:K}
  \begin{tabular}{c||cc|cc||cc|cc}
    &\multicolumn{4}{c||}{$\kc\in\PS{0}(\Oh)$}
    &\multicolumn{4}{c}  {$\kc\in\PS{1}(\Oh)$}\\
    \hline
    \textsf{cells}
    &\multicolumn{2}{c|}{\SchemeI}&\multicolumn{2}{c||}{\SchemeII}
    &\multicolumn{2}{c|}{\SchemeI}&\multicolumn{2}{c}  {\SchemeII}\\
    \hline
    {} 
    &\textsf{err(p)} & \textsf{err(\k\bu)} & \textsf{err(p)} & \textsf{err(\k\bu)}
    &\textsf{err(p)} & \textsf{err(\k\bu)} & \textsf{err(p)} & \textsf{err(\k\bu)}\\
    \hline
    412     & 2.762e-3 & 7.451e-3 & 5.438e-3 & 2.679e-2 & 2.903e-3 & 3.063e-3 & 2.588e-3 & 3.067e-3\\
    1591    & 6.976e-4 & 2.370e-3 & 3.271e-3 & 1.426e-2 & 6.540e-4 & 9.656e-4 & 6.541e-4 & 9.637e-4\\
    6433    & 1.650e-4 & 9.264e-4 & 2.242e-3 & 7.354e-3 & 1.548e-4 & 4.897e-4 & 1.548e-4 & 4.887e-4\\
    25698   & 4.066e-5 & 2.581e-4 & 1.104e-3 & 3.690e-3 & 3.833e-5 & 1.267e-4 & 3.832e-5 & 1.267e-4\\
    102772  & 1.007e-5 & 1.134e-4 & 4.802e-4 & 2.076e-3 & 9.502e-6 & 5.545e-5 & 9.502e-6 & 5.543e-5\\
    \hline
    \textsf{rate} & 2.04 & 1.53 & 0.86 & 0.94 & 2.06 & 1.45 & 2.03 & 1.45
  \end{tabular}
\end{table}


\section{Conclusions}
\label{sec:conclusions}

Numerical schemes for nonlinear parabolic equations based on harmonic
averaging of cell-centered diffusion coefficients at cell interfaces
break down when some of these coefficients go to zero or their ratio
is too large.
To address this issue,
in~\cite{Lipnikov-Manzini-Moulton-Shashkov:2016} we proposed a new
family of second-order accurate mimetic finite difference schemes on
polygonal and polyhedral meshes.
In this new discrete setting the primary mimetic operator approximates
the continuum operator $\textsf{div}(\k\,\cdot\,)$, while the derived
(dual) mimetic operator approximates $\nabla(\cdot)$.
The discrete divergence operator requires a staggered discretization
of the diffusion coefficient, one value per mesh cell and up to two
values per mesh face.
The availability of face diffusion coefficients provides more
flexibility in the numerical formulation, which can be exploited to
design robust numerical algorithms for nonlinear problems.
For instance, upwinding of the diffusion coefficients on mesh faces
can be easily incorporated into the new mimetic schemes.
The new mimetic method applied to the steady diffusion equation in
mixed form is proved to be well-posed since it satisfies the discrete
inf-sup condition, and convergent by deriving first-order error
estimates for the scalar and gradient unknowns.
Numerical experiments verify the theory.




\vskip3mm
\section*{Acknowledgments}

This work was carried out under the auspices of the
National Nuclear Security Administration of the U.S. Department of
Energy at Los Alamos National Laboratory under Contract
No. DE-AC52-06NA25396.
The authors acknowledge the support of the US Department of Energy
Office of Science Advanced Scientific Computing Research (ASCR)
Program in Applied Mathematics Research.
The article is assigned the LA-UR number LA-UR-16-28012.
We are grateful to Dr. Rao Garimella (LANL) for helping us with
generating dual meshes satisfying geometric constraints.
All meshes in this paper were created using his mesh generation
toolset MSTK (software.lanl.gov/MeshTools/trac).


\vskip4mm
\bibliography{theory_SINUM}
\bibliographystyle{plain}

\end{document}